\documentclass[12pt,a4paper]{article}%
\usepackage{amsmath}
\usepackage{graphicx}
\usepackage{epsfig}%
\usepackage{amsfonts}%
\usepackage{amssymb}
\usepackage{algpseudocode}						
\usepackage{algorithm}	
\usepackage{color}
\usepackage{mathtools}
\usepackage{diagbox}						

\newtheorem{theorem}{Theorem}[section]

\newtheorem{lemma}[theorem]{Lemma}

\newtheorem{remark}[theorem]{Remark}

\makeatletter\@addtoreset{equation}{section}\makeatother
\makeatletter\@addtoreset{figure}{section}\makeatother
\makeatletter\@addtoreset{table}{section}\makeatother
\newcommand{\dotprod}[2]{\langle#1,#2\rangle} 					
\newcommand{\R}{{\mathbb{R}}}       							

\DeclareMathOperator{\diag}{diag}								

\begin{document}

\title{Tensor product method for fast solution of optimal control problems with 
fractional multidimensional Laplacian in constraints}
 
\author{Gennadij Heidel\thanks{University of Trier,
        FB 4 - Department of Mathematics, D-54286, Trier,
        Germany ({\tt   shaw@uni-trier.de}).}
        \and
        Venera Khoromskaia\thanks{Max Planck Institute for
        Mathematics in the Sciences, Leipzig, Germany ({\tt vekh@mis.mpg.de}).}
        \and
        Boris N. Khoromskij\thanks{Max-Planck-Institute for
        Mathematics in the Sciences, Inselstr.~22-26, D-04103 Leipzig,
        Germany ({\tt bokh@mis.mpg.de}). 
        }
        \and
        Volker Schulz\thanks{University of Trier,
        FB 4 - Department of Mathematics, D-54286, Trier,
        Germany ({\tt Volker.Schulz@uni-trier.de}).}
        }

 \date{}

\maketitle

\begin{abstract}
 We introduce the tensor numerical method for  solution of the $d$-dimensional optimal 
 control problems with fractional Laplacian type 
operators in constraints discretized on large $n^{\otimes d}$ tensor-product Cartesian grids.
The approach is based on the rank-structured approximation 
of the matrix valued functions of the corresponding fractional finite difference Laplacian.  
We solve the equation for the control function, where the system matrix includes the sum of 
the fractional $d$-dimensional Laplacian and its inverse. The matrix valued functions of 
discrete  Laplace operator on a tensor grid are  diagonalized by using the fast Fourier 
transform (FFT). Then the low rank   approximation of the $d$-dimensional tensors obtained by 
folding of the corresponding large diagonal matrices of eigenvalues are computed,
which allows to solve the governing equation for the control function in a tensor-structured format.  
The existence of low rank canonical approximation to the class of 
matrix valued functions involved is justified  by using the
sinc quadrature approximation method applied to the Laplace transform of the generating function.
The linear system of equations for the control function is solved by the PCG iterative 
method with the rank truncation at each iteration step, where the low Kronecker rank preconditioner 
is precomputed. The right-hand side, the solution vector, and the governing system matrix are maintained 
in the rank-structured tensor format which beneficially reduces the numerical cost to $O(n\log n)$, 
outperforming  the standard FFT based methods of complexity $O(n^3\log n)$ for 3D case.  
Numerical tests for the 2D and 3D control problems confirm the linear complexity scaling of the 
method in the univariate grid size $n$.  
\end{abstract}

\section{Introduction}\label{Int:SH}

 Optimization problems with partial differential equations (PDEs) as constraints 
 are well known in the mathematical literature and have been studied for 
many years, see \cite{Pontryag:1963,Troeltzsch:2005,Hinze_book:09} 
for comprehensive presentation.
 A major challenge in the numerical analysis  of the PDE constrained optimal control problems
is the evaluation of the constraints requiring 
the solution of a partial differential or nonlocal integral-differential equations in $\mathbb{R}^d$. 
Therefore, to make these problems tractable, specially tailored solvers are required.
In the classical sense, partial differential equations are governed  by local operators, 
implying the sparsity of the discrete operator (matrix). In this case,
traditional elliptic problem solvers for the forward problem can be modified 
to apply to the optimization problems, see \cite{Benner:04,HerKun:2010} for an overview. 
Multigrid methods for elliptic equations are shown to be efficient since their 
computational complexity is linear in the number of grid points in the computational 
domain in $\mathbb{R}^d$, see \cite{BoSch:2009,BoSch:2012}.

Recently,  study of the  problems with nonlocal constraints, where the operator is  
of integral type, has attracted growing interest. In particular, the problems  
 including   the \emph{fractional Laplacian} operator $(-\Delta)^{\alpha}$, 
for $0< \alpha \leq 1$, can be applied for modeling of 
complex systems \cite{Ainsworth:18,Bonito:18}. 
 Several definitions of fractional elliptic operator based on either integral 
or spectral formulations have been considered in the literature
\cite{Podlub:99,GaHaKh3:02,HaHighTrefeth:08,Kwasn:17,FracLapl:2018,RapJel:10}. 
 The main distinction between the integral and spectral definitions is that
 the first one admits a nonlocal boundary condition, 
while the spectral formulation requires only the local boundary condition \cite{FracLapl:2018}.
 A number of application fields motivating the use of fractional power of elliptic operators,
 for example in biophysics, mechanics, nonlocal electrostatics and image processing 
 have been discussed
in \cite{AtPSZo:14,AntilOtarola,DPSS:2014,Vabich:2015,DuLazPas:18,HLMMV:2018,KaMelenk:19}. 
 In such applications control problems arise consistently.

 In order to illustrate our problem formulation, we notice that in \cite{AntilOtarola}
 it is shown that heat diffusion in special materials, 
 called "plasmonic nanostructure networks" \cite{Ben_Abd_Messina:13}, 
 which have strong relation to new composite 
 materials, is described exactly by fractional Laplacian operator. 
 This means that the interpretation for the case $\alpha=1$ in the form of  
 distributed control of the heat equation can be carried immediately over 
 to the control of heat equation in the particular material (see also \cite{RapJel:10}).
 In our particular application one of the motivations to use the fractional elliptic operators 
${\cal L}^\alpha$, $\alpha\in (0,1]$, in constraints  
 is due to the opportunity to control the sharpness (contrast) for the resolution  
 of the control and design functions by tuning the fraction power $\alpha$ in  a direction to a 
 smaller value.

Numerical solution of problems with nonlocal operators poses an additional
computational challenge: since local information is not sufficient for the evaluation of the operator,
the discretized operator will be a dense matrix instead of a sparse one--if implemented in a
straightforward way.
A number of approaches have been proposed and analyzed in the literature 
to circumvent these difficulties, see papers 
\cite{Higham_MatrFunc:08,HaHighTrefeth:08,Vabich:2015,DuLazPas:18,HLMMV:2018,Schwab:18}
considering approximation methods for fractional elliptic operators,
and \cite{DPSS:2014,Vabi:18} related to time dependent problems.  
An extension method   proposed in \cite{CaffSilv}, 
reduces a fractional Laplacian problem to a classical
Laplacian problem in a higher-dimensional space, and allows to  make the 
problem tractable in some cases, see also recent paper \cite{Schwab:18}, 
which develops the higher-order approximation methods. 
A proof of concept for an optimal control solver
based on the extension approach is considered in \cite{AntilOtarola}.

 The above methods   based on the standard numerical techniques   
provide at least linear complexity scaling in the discrete problem size,  
thus exhibiting an exponential increase   of storage and computational complexity  
in the number of dimensions $d$, as $O(n^d)$, where $n$ is the univariate grid size 
 for discretizations over product $n^{\otimes d}$ grids in $\mathbb{R}^d$.  
Modern tensor numerical methods allow to reduce the computational cost for 
multivariate problems from $O(n^d)$ to $O(n)$. 
They appeared as a bridging of the tensor decompositions 
from multilinear algebra  \cite{DMV-SIAM2:00,Cichocki:2002,Hack_Book12} 
with basic results in nonlinear approximation theory  on  
the low-rank separable 
representation of the multivariate functions and operators \cite{GHK:05,GaHaKh3:02}.   
The main idea of  tensor  numerical methods is to reduce the solution procedures for  
the multivariate nonlocal integral-differential equations to tensor product 
operations on one-dimensional data.  
 Nowadays, the development and application of low-rank tensor numerical methods  is one of the 
prior  directions  in scientific computing 
\cite{KhSch:2011,KreTob:2010,DolOs:2011,Cichocki:2016,KhorBook:17,Khor2Book2:18,MaRaSch:19}.

In this article\footnote{The present paper is an extended version of the preprint \cite{HKKS:18}.}, 
we introduce the tensor numerical approach for the efficient solution 
of the  optimal control problems with fractional $d$-dimensional Laplacian type 
operators in constraints. 
We solve the equation for the control function, where the system matrix includes the sum of 
the fractional $d$-dimensional Laplacian and its inverse.
We propose and analyze the approximate low-rank structure representations of functions 
of the fractional Laplacian $(-\Delta)^{\alpha}$, and its inverse $(-\Delta)^{-\alpha}$, $\alpha >0$, 
in the bounded domain in $\mathbb{R}^d$, $d=2,3$, by using the canonical tensor format.
 Algorithmically, the functions of discretized  Laplace operator are diagonalized
by the FFT transform, 
and then the low rank   approximation 
 of the $d$-dimensional tensors obtained by folding of the corresponding large diagonal matrices 
 of eigenvalues are computed.
This  allows us to solve the governing equation for the control function in a tensor-structured format.
Given the low Kronecker rank preconditioner, 
the PCG iteration with rank truncation by using the representation of the operator system matrix  
  in the canonical tensor format is introduced.
In the three-dimensional case, the multigrid Tucker tensor decomposition  is applied, 
  while the canonical rank reduction in a course of iteration
is performed by using the tensor transforms based on the Reduced Higher 
Order SVD  (RHOSVD) \cite{KhKh3:08,Khor2Book2:18}.

The theoretical justification for the low-rank canonical and Tucker tensor approximation 
of functions of the discrete fractional Laplacian in $\mathbb{R}^d$, $d=2,3, $ is provided. 
The existence of the low-rank canonical decompositions is based on the theory of 
the sinc-approximation methods applied to the Laplace transform representation of the respective
operator/matrix values functions. 
We show that these low-parametric representations transfer to the solution operator 
for equation describing the control variable, which includes the sum of a fractional Laplacian
and its inverse.  
In this way the spacial dimensions can be approximately separated resulting in a low-rank 
tensor structure in the solution vector,  provided that the design function in 
the right-hand side is approximated with the low rank. 
 
 The tensor approach introduced in this paper  is based on the spectral decomposition of the target operators
  and applies to the problems discretized on large $n\times n \times n$ Cartesian grids 
in the 3D box-type domains. The low rank tensor structure of 
all matrix valued functions involved allows to reduce the numerical cost to the linear scale 
$O(n \log n)$ in the univariate problem size, contrary to the traditional numerical schemes 
based on FFT diagonalization, which amount to
the linear-logarithmic  complexity in the volume size, $O(n^3 \log n)$. 
Numerical tests for the 2D and 3D control problems 
confirm the linear complexity scaling of the method in the univariate grid size $n$. 
Moreover, 
the efficient representation of matrix valued functions of fractional Laplacian can be used
in various applications as the spectrally close preconditioner for the control problems 
with the more general elliptic operators in the constraint, and for more general geometries.  

 Notice that the low-rank tensor method  for solution of a fractional time dependent 
optimal control problems, where the operator in constraints is constructed 
as a sum of  the Riemann-Lioville fractional one-dimensional Laplacian operators 
(integral Caputo type representation discretized by using one-level Toeplitz matrices)  
has been considered in \cite{DPSS:2014}.
See Remarks \ref{rem:direct_fract_Lap} and \ref{rem:GenerLapl} for the discussion on the special case of ``directionally fractional'' Laplacian in constraints.

 The rest of the paper is organized as follows. Section \ref{Int2:SH} describes the 
 considered problem setting. Section \ref{sec:Larrange_eqn} discusses the FEM/FDM 
 discretization schemes for functions of an elliptic operator,
 formulates the traditional Lagrange multiplies approach and describes the Kronecker product 
 tensor structure in the functions of discrete Laplacain in $d$ dimensions.
  Section \ref{sec:tensor_approx_theory} analyzes the tensor approximation of the inverse 
 to fractional Laplace operator and 
  some other matrix valued functions of fractional Laplacian in $\mathbb{R}^d$ 
 arising in representation of the unknown control and  design functions.
 In Section \ref{sec:low-rank_tensor} we recall the main definition and basic properties of the canonical 
 and Tucker tensor formats to be applied for tensor approximation and rank truncation.
 Finally, in Section \ref{sec:numerics_tensor}, we collect the results of numerical 
 tests for 2D and 3D examples which confirm the efficiency of the tensor approximation for the 
 considered class of optimal control problems.
 
 

\section{Problem setting }\label{Int2:SH} 

Our goal is the construction of fast solution schemes for solving the control problems
with $d$-dimensional fractional elliptic operators in the constraints. 
For this reason we confine ourself to the case of box-type domains and to the class
of generalized Laplacian type elliptic operators with separable coefficients.


Given the design function $y_\Omega \in L^2(\Omega)$ on $\Omega :=(0,1)^d$, $d=1,2,3$,
first, we consider the optimization problem for the cost functional
\begin{equation} \label{eqn:cost_func}
 \min_{y,u} J(y,u):=\int_\Omega (y(x) -y_\Omega(x))^2\, dx + \frac{\gamma}{2} 
 \int_\Omega u^2(x) \,dx, \quad \gamma>0, 
\end{equation}
constrained by the elliptic boundary value problem  
in $\Omega$ for the state variable  $y \in H_0^1(\Omega)$,  and control variable $u$, 
\begin{equation} \label{eqn:basic_setting}
  {\cal A} y  := -\nabla^T \cdot \mathbb{A}(x)\nabla y = \beta u, \quad x\in \Omega,\;
  u\in L_2(\Omega),\quad \beta >0,
  \end{equation}
endorsed with the homogeneous Dirichlet boundary conditions on $\Gamma = \partial \Omega $, 
i.e., $y_{|\Gamma}=0$.
The coefficient matrix $\mathbb{A}(x)\in \mathbb{R}^{d\times d}$ is supposed to 
be symmetric, positive definite and uniformly bounded in $\Omega $ 
with positive constants $c>0$ and $C>0$, i.e.,
\[
 c\, I_{d\times d}   \leq \mathbb{A}(x)\leq C\, I_{d\times d}.
\]

Under above assumptions the associated bilinear form
\[
 A(u,v)=\int_\Omega \mathbb{A}(x)\nabla u(x) \cdot \nabla v(x) \, dx
\]
defined on $V\times V $, $V:=\{v\in H_0^1(\Omega) \}$ is symmetric, coercive and bounded on $V$
with the same constants $c$ and $C$.

In what follows, we describe the tensor method for fast numerical solution of the 
optimization problem with the generalized constraints
\begin{equation} \label{eqn:frac_setting}
  {\cal A}^\alpha y  = \beta u(x), \quad x\in \Omega,
\end{equation}
such that for $0 < \alpha \leq 1$, the fractional power of the elliptic operator ${\cal A}$ is defined 
by
\begin{equation} \label{eqn:frac_Oper}
 {\cal A}^\alpha y(x) = \sum_{i=1}^\infty \lambda_i^\alpha c_i \psi_i(x), \quad 
 y=\sum_{i=1}^\infty c_i \psi_i(x),
\end{equation}
where $\{\psi_i(x)\}_{i=1}^\infty$ is the set of $L_2$-orthogonal eigenfunctions
of the symmetric, positive definite operator ${\cal A}$, while $\{\lambda_i\}_{i=1}^\infty$ are 
the corresponding (real and positive) eigenvalues.

In the present paper, we consider the particular case of fractional Laplacian in the form 
\begin{equation}  \label{eqn:frac_setting_Lap}
{\cal A}^\alpha := 
 (\sum\limits_{\ell=1}^d -\Delta_\ell )^\alpha, \quad \alpha>0,
\end{equation} 
where $\Delta_\ell$ is the 1D Laplacian in variable $x_\ell$. 

Notice that the elliptic operator inverse ${\cal A}^{-1}={\cal T}:L_2(\Omega) \to V$, where 
${\cal A}={\cal T}^{-1}$,
provides the explicit representation for the state variable, 
$y=\beta {\cal T} u = \beta {\cal A}^{-1} u$ in case 
(\ref{eqn:basic_setting}), while in the general case (\ref{eqn:frac_setting}) we have
\begin{equation} \label{eqn:control}
y =\beta {\cal T}^{\alpha} u \equiv \beta {\cal A}^{-\alpha} u.
 \end{equation}
 Here ${\cal T}$ is a compact, symmetric and
positive definite operator on $L_2(\Omega)$ and its eigenpairs $\{\psi_i,\mu_i\}$, 
$i=1,\dots,\infty$, provide an orthonormal basis for $L_2(\Omega)$, where we have 
$\mu_i = \lambda_i^{-1}$.

There are several representations (definitions) for the fractional power 
of the symmetric, positive definite operators
${\cal A}^{\alpha}$ and ${\cal A}^{-\alpha}$ with $0< \alpha \leq 1$, see for example
the survey  papers \cite{Kwasn:17,FracLapl:2018}. 
In particular, the Dunford-Taylor-Cauchy contour integral, the Laplace transform 
integral representations  and the spectral definitions could be applied.

In the presented computational schemes based on low rank tensor decompositions,
 we apply the spectral definition and use the Laplace transform integral representation
for the analysis and justification of the low rank tensor approximation. 
For $\alpha >0$,
the integral representation based on the Laplace transform
\begin{equation} \label{eqn:Laplace_transf}
 {\cal A}^{-\alpha}=\frac{1}{\Gamma(\alpha)} \int_0^\infty t^{\alpha -1} e^{-t {\cal A}}\, dt
\end{equation}
suggests the numerical schemes for low rank canonical tensor representation of the 
operator (matrix) ${\cal A}^{-\alpha}$ by using the
sinc quadrature approximations for the integral on the real axis \cite{GHK:05}, 
\begin{equation} \label{eqn:sinc_Laplace_tr}
  \int_0^\infty t^{\alpha -1} e^{-t {\cal A}}\, dt \approx 
 \sum\limits_{k=-M}^M \hat{c}_k t_k^{\alpha -1} e^{-t_k {\cal A}}= \sum\limits_{k=-M}^M c_k 
\bigotimes_{\ell=1}^d e^{-t_k {A_\ell}},
\end{equation}
applied to the operators composed by a sum of commutable terms, 
$$
{\cal A}=\sum_{\ell=1}^d A_\ell, \quad [A_\ell,A_k]=0,\quad \mbox{for all} \quad \ell,k=1,\cdots,d,
$$
which ensures that each summand in (\ref{eqn:sinc_Laplace_tr}) is separable, i.e.
$e^{-t_k {\cal A}}=\bigotimes_{\ell=1}^d e^{-t_k {A_\ell}}$.
For example, in the case of Laplacian we have the $d$-term decomposition
$\Delta = \sum_{\ell=1}^d \frac{\partial^2}{\partial x_\ell^2}$ with commutable 1D operators. 
The efficiency of this tensor approach is justified by the established exponentially fast 
convergence of the sinc quadratures in the number of separable terms $M$ in (\ref{eqn:sinc_Laplace_tr}), 
applied on a class of analytic functions. 
  In the present paper,  this techniques is used for both the theoretical 
analysis of the rank decomposition schemes and for the description of their
constructive representation,  see Section \ref{sec:tensor_approx_theory}.

Further more, for $\alpha >0$ the Dunford-Taylor (or Dunford-Cauchy) contour integral
representation reads (see for example \cite{Higham_MatrFunc:08,HaHighTrefeth:08,GaHaKh3:02})
\begin{equation} \label{eqn:Dunford_Int}
 {\cal A}^{-\alpha}=\frac{1}{2\pi i} \int_{\cal G} z^{-\alpha} ({\cal A} -z{\cal I})^{-1}dz,
\end{equation}
where the contour ${\cal G}$ in the complex plane 
includes the spectrum of operator (matrix) ${\cal A}$.
This representation applies to any $u\in L_2(\Omega)$ and it allows to define the negative 
fractional powers of elliptic operator as a finite sum of elliptic resolvents 
${\cal R}_z({\cal L})= ( z{\cal I} - {\cal L})^{-1}$
by application of certain quadrature approximations,  
\[
 \int_{\cal G} z^{-\alpha} ({\cal A} -z{\cal I})^{-1}dz 
 \approx  \sum\limits_{k=1}^M c_k z_k^{-\alpha}({\cal A} -z_k{\cal I})^{-1}, \quad z_k \in {\cal G},
\]
similar to (\ref{eqn:sinc_Laplace_tr}),  see also \cite{HaHighTrefeth:08,GHK:05,GaHaKh3:02}. 
This opens the way for 
multigrid based, ${\cal H}$-matrix (see \cite{GaHaKh3:02}) or tensor-structured schemes approximating the 
fractional power of elliptic operator with variable coefficients of rather general form
  and defined on non-rectangular geometries  combined with the nested iterations.

It is worth to notice that both integral representations (\ref{eqn:Laplace_transf}) and
(\ref{eqn:Dunford_Int}) can be applied to rather general class of 
analytic functions of operator $f({\cal A})$, including the case $f({\cal A})=\mathcal{A}^{-\alpha}$, 
see \cite{Higham_MatrFunc:08,GHK:05,GaHaKh3:02}.

The constraints equation (\ref{eqn:control}) allows to derive the Lagrange 
equation for the control $u$ in the explicit form as follows 
(see \S\ref{sec:Larrange_eqn} concerning the Lagrange equations)
\begin{equation} \label{eqn:Lagrange_cont}
 \big(\beta\mathcal{A}^{-\alpha} + \tfrac{\gamma}{\beta}\mathcal{A}^{\alpha}\big) u = y_\varOmega,
\end{equation}
for some positive constants $\beta>0$ and $\gamma>0$. This equation implies the following 
representation for the state variable
\begin{equation} \label{eqn:State}
y= \beta\mathcal{A}^{-\alpha} u. 
\end{equation}
The practically interesting range of parameters includes the case $\beta=O(1)$ for the small 
values of $\gamma>0$.
Our  tensor numerical method  is designed for  solving 
equations (\ref{eqn:Lagrange_cont}) and (\ref{eqn:State}) that include the nonlocal operators
of ``integral-differential'' type. 
The efficiency of the rank-structured tensor approximations presupposes 
that the design function in the right-hand side of these equations,
$y_\varOmega(x_1,x_2,x_3)$, allows the low rank separable approximation.

Since we aim for the low-rank (approximate) tensor representation of all functions and 
operators involved in the above control problem, 
it is natural to assume that the equation coefficients matrix takes a diagonal form
\[
 \mathbb{A}(x) = \mbox{diag}\{a_1(x_1),a_2(x_2) \}, \quad a_\ell(x_\ell) >0, 
 \quad \ell=1,2,
\]
in 2D case, and similar for $d=3$
\begin{equation} \label{eqn:DiagCoef3}
 \mathbb{A}(x) = \mbox{diag}\{a_1(x_1),a_2(x_2),a_3(x_3) \}, \quad a_\ell(x_\ell) >0, 
 \quad \ell=1,2,3,
\end{equation}
where the $d$-dimensional Laplacian is the particular case with $a_\ell(x_\ell)=const$.

  In what follows, we consider the discrete matrix formulation of the  optimal control problem 
  (\ref{eqn:cost_func}),  (\ref{eqn:frac_setting})
 based on the FEM/FDM discretization ${A}_h$ of $d$-dimensional Laplacian  
 defined on the uniform $n_1\times n_2 \times \ldots \times n_d  $ tensor grid in $\Omega$,
  where $h=h_\ell=1/n_\ell$ is the univariate mesh parameter. 
  The $L_2$ scalar product 
  is substituted by the Euclidean scalar product $(\cdot,\cdot)$ of vectors in 
  $\mathbb{R}^{\bf n}$, ${\bf n}=(n_1, n_2,\ldots, n_d)$. 
  
  In this paper, the fractional Laplacian ${\cal A}^{\alpha}$ 
  is approximated by its FEM/FDM representation $({\cal A}_h)^{\alpha}$, where
  the matrix $({\cal A}_h)^{\alpha}$ is defined by spectral decomposition 
  of  ${A}_h$ in (\ref{eqn:frac_Oper}). 
  The FEM approximation theory for fractional powers of elliptic operator was presented in 
 \cite{Ainsworth:18,DuLazPas:18,HLMMV:2018,KaMelenk:19}, see also literature therein.

 \section{Optimality conditions and representations in a low-rank format}
\label{sec:Larrange_eqn}

The solution of problem \eqref{eqn:cost_func} with constraint \eqref{eqn:frac_setting} requires
solving for the necessary optimality conditions. In this section, we will derive these conditions
based on a discretize-then-optimize-approach. Then, we will discuss how the involved discretized
operators can be applied efficiently in a low-rank format, and how this can be used to design
a preconditioned conjugate gradient (PCG) scheme for the necessary optimality conditions.

\subsection{Discrete optimality conditions}
We consider the discretized version of the control problem 
\eqref{eqn:cost_func}-\eqref{eqn:frac_setting}. We assume we have a uniform grid in each space dimension,
i.\,e. we have $N=n_1 n_2$ (for $d=2$) or $N=n_1 n_2 n_3$ (for $d=3$) grid points. We
will denote the discretized state $y$, design $y_\Omega$ and control $u$ by vectors
$\mathbf{y}, \mathbf{y}_\Omega, \mathbf{u}\in \R^N$, respectively. For simplicity, 
we assume that we use the same approximation for all quantities. 

Then, the discrete problem is given as

\begin{align*}
	\min_{\mathbf{y},\mathbf{u}} =  &\frac{1}{2}(\mathbf{y} - \mathbf{y}_\Omega)^TM(\mathbf{y} 
	- \mathbf{y}_\Omega)  + \frac{\gamma}{2} \mathbf{u}^T M \mathbf{u}\\
	\text{s.\,t.} ~ A^\alpha \mathbf{y} = & \beta M\mathbf{u},
\end{align*}
where $A = \mathcal{A}_h$ denotes a discretization of the elliptic operator $\mathcal{A}$ by 
finite elements or finite differences. 
The matrix $M$ will be a mass matrix in the finite element case and simply the 
identity matrix in the finite difference case.

For the discrete adjoint $\mathbf{p}$ define the Lagrangian function
\begin{equation}
	L(\mathbf{y},\mathbf{u},\mathbf{p}) = \frac{1}{2}(\mathbf{y} - 
	\mathbf{y}_\Omega)^T M(\mathbf{y} - \mathbf{y}_\Omega)  
	+ \frac{\gamma}{2} \mathbf{u}^T M \mathbf{u} + \mathbf{p}^T(A^\alpha \mathbf{y} - 
	\beta M \mathbf{u}),
\end{equation}
and compute the necessary first order conditions, given by the Karush-Kuhn-Tucker (KKT) system,

\begin{equation}\label{eqn:KKT_for_yup}
	\begin{bmatrix}
		M & O &  A^\alpha\\
		O & \gamma M & -\beta M\\
		A^\alpha & -\beta M & O
	\end{bmatrix}
	\begin{pmatrix}
		\mathbf{y}\\
		\mathbf{u}\\
		\mathbf{p}\\
	\end{pmatrix} =
	\begin{pmatrix}
		\mathbf{y}_\varOmega\\
		\mathbf{0}\\
		\mathbf{0}\\
	\end{pmatrix}.
\end{equation}

We can solve the state equation for $\mathbf{y}$, getting
\begin{equation*}
	\mathbf{y} = \beta A ^{-\alpha}\mathbf{u},
\end{equation*}
and the design equation for $\mathbf{p}$ getting
\begin{equation*}
	\mathbf{p} = \tfrac{\gamma}{\beta} \mathbf{u} .
\end{equation*}
Hence the adjoint equation gives us an equation for the control $\mathbf{u}$, namely
\begin{equation}\label{eqn:forControl}
A_2 \mathbf{u}\equiv 
\big(\beta A^{-\alpha} + \tfrac{\gamma}{\beta}A^\alpha\big) \mathbf{u}  = \mathbf{y}_\varOmega.
\end{equation}

In this paper, we use the explicit relation (\ref{eqn:forControl})  as the  
governing equation for calculation of the unknown control $\mathbf{u}$.
The main motivation is that the system matrix in equation (\ref{eqn:forControl}),
that includes fractional power of the stiffness matrix $A$ and its inverse $A^{-1}$, 
can be well approximated by the low Kronecker rank matrices. 
Hence, the reformulation of the traditional system (\ref{eqn:KKT_for_yup}) in the form
(\ref{eqn:forControl}) proves to be convenient for the construction of low-rank tensor 
approximation to the arising hybrid system matrix  $A_2$ and for the construction of 
spectrally close preconditioner  by approximation of $A_2^{-1}$. Furthermore,
notice that the proposed  rank-structured approximations for the case of Laplacain type 
operators can be also gainfully used for preconditioning in the case of more general elliptic
operators with variable coefficients in constraints.
In what follows, we construct and analyze the low Kronecker 
rank decompositions of the arising matrix valued functions, see \S \ref{ssec:Multivariate_cores} 
and \S \ref{ssec:Theory_multivariate_cores} below. 

\subsection{Matrix-vector multiplication in the low-rank format}

 First, we discuss a Kronecker form decomposition of functions of discrete Laplacian $A=A_h$ 
diagonalized in the Fourier basis, 
which is compatible with low-rank data. Let $I_{\ell}$ denote the identity matrix, and $A_{(\ell)}$
the discretized one-dimensional Laplacian on the given grid in the $\ell$-th mode, then we have
\begin{equation} \label{eqn_low_rank_lap}
	A = A_{(1)} \otimes I_2 \otimes I_3 + I_1 \otimes A_{(2)} \otimes I_3 + 
	I_1 \otimes I_2 \otimes A_{(3)},
\end{equation}
where $\otimes$ denotes the Kronecker product of matrices.
To calculate the matrices $A_{(\ell)}$, we simply discretize the one-dimensional subproblems
\begin{align*}
	-&y''(x_\ell) = \beta u (x_\ell)\\
	&y(0) = 0 = y(1).
\end{align*}
Using a uniform grid with grid size $h_\ell$, we obtain the discretizations
\begin{equation*}
\underbrace{-\frac{1}{h_\ell}
\begin{bmatrix}
	2   & -1      &       &       \\
	-1   & \ddots     & \ddots      &  \\
    & \ddots & \ddots & -1       \\
    &        &      -1  & 2      
\end{bmatrix}}
_{=A_{(\ell)}}
\begin{pmatrix}
	y^{(\ell)}_1\\
	\vdots\\
	y^{(\ell)}_{n_\ell}
\end{pmatrix}
=\beta
\begin{pmatrix}
	u^{(\ell)}_1\\
	\vdots\\
	u^{(\ell)}_{n_\ell}
\end{pmatrix}
\end{equation*}

The ``one-dimensional'' matrix $A_{(\ell)}$ has an eigenvalue decomposition in the Fourier basis, i.e.
\begin{equation*}
	A_{(\ell)} = F_\ell^* {\varLambda}_{(\ell)} F_\ell.
\end{equation*}
In the case of homogeneous Dirichlet boundary conditions the matrix $F_\ell$ defines the 
$\sin$-Fourier transform while 
${\varLambda}_{(\ell)}=\mbox{diag}\{\lambda^{(\ell)}_1,\dots,\lambda^{(\ell)}_n\}$, where 
$\lambda_k$ denote the eigenvalues of the univariate discrete Laplacian with 
Dirichlet boundary conditions. These are given by
\begin{equation}\label{eqn:eig_lap}
	\lambda_k = -\frac{4}{h_\ell^2}\sin^2 \bigg( \frac{\pi k}{2(n_\ell+1)} \bigg)
	 = -\frac{4}{h_\ell^2}\sin^2 \bigg( \frac{\pi k h_\ell}{2} \bigg).
\end{equation}

Thus, using the properties of the Kronecker product, we can write the first 
summand in \eqref{eqn_low_rank_lap} as
\begin{equation*}
\begin{split}
A_1 \otimes I_2 \otimes I_3 &= (F_1^* {\varLambda}_{(1)} F_1) \otimes (F_2^* I_2 F_2) \otimes (F_3^* I_3 F_3) \\
&= (F_1^*\otimes F_2^* \otimes F_3^*) ({\varLambda}_{(1)}\otimes I_2\otimes I_3) (F_1\otimes F_2 \otimes F_3).
\end{split}
\end{equation*}
The decomposition of the second and third summand works analogously, thus we can write equation 
\eqref{eqn_low_rank_lap} as
\begin{equation}\label{eqn:DiagLaplace}
\begin{split}
A =  &(F_1^*\otimes F_2^*\otimes F_3^*) (\varLambda_1\otimes I_2\otimes I_3) (F_1\otimes F_2 \otimes F_3)\\
	&+(F_1^*\otimes F_2^*\otimes F_3^*) (I_1\otimes \varLambda_2\otimes I_3) (F_1\otimes F_2 \otimes F_3)\\
	&+(F_1^*\otimes F_2^*\otimes F_3^*) (I_1\otimes I_2\otimes \varLambda_3) (F_1\otimes F_2 \otimes F_3)\\
	=\; &(F_1^*\otimes F_2^*\otimes F_3^*) 
	\underbrace{\big( \varLambda_1\otimes I_2\otimes I_3 + I_1\otimes \varLambda_2\otimes I_3 + 
	I_1\otimes I_2\otimes \varLambda_3 \big)}_{\eqqcolon \varLambda}
	(F_1\otimes F_2 \otimes F_3),
\end{split}
\end{equation}
 where $\varLambda \in \mathbb{R}^{n^3 \times n^3} $ is the diagonal matrix.
The above expression gives us  the eigenvalue factorization (diagonalization) of an
analytic matrix valued functions ${\cal F}(A)$ of the matrix $A$,
\begin{equation}\label{eqn:DiagFLaplace3D}
\mathcal{F}(A)  = (F_1^*\otimes F_2^*\otimes F_3^*) 
	\mathcal{F}(\varLambda) (F_1\otimes F_2 \otimes F_3),
\end{equation}
which can be multiplied with a vector of general structure in $O(n^3 \log n)$ operations.

In the tensor based scheme,
we construct the low Kronecker rank tensor decomposition of the diagonal matrix 
$\mathcal{F}(\varLambda)\in \mathbb{R}^{n^3 \times n^3}$,
\begin{equation*}
	\mathcal{F}(\varLambda)  = \sum_{k=1}^R 
	\diag \big(\mathbf{u}_1^{(k)} \otimes \mathbf{u}_2^{(k)} \otimes \mathbf{u}_3^{(k)} \big),
\end{equation*}
with vectors $\mathbf{u}_\ell^{(k)}\in\R^{n_\ell}$ and $R \ll \min(n_1,n_2, n_3)$,
via canonical approximation of the 
$n\times n \times n$ tensor obtained by folding of the long vector 
$\mbox{diag}(\mathcal{F}(\varLambda))\in \mathbb{R}^{n^3}$.
Then, the low Kronecker rank approximation for $\mathcal{F}(A)$ 
is obtained by multiplication from left and right with the factorized FFT matrices,
see Section \ref{sec:tensor_approx_theory}. 


 First, we consider the fast calculation of the matrix vector product with $A$. 
In the case $d=2$, the factorization (\ref{eqn:DiagLaplace}) simplifies to
\begin{equation}\label{eqn:DiagLaplace2D}
\begin{split}
A       &= (F_1^*\otimes F_2^*) 
	\underbrace{\big( \varLambda_1\otimes I_2 + 
	I_1\otimes \varLambda_2 \big)}_{\eqqcolon \varLambda} 	(F_1\otimes F_2 ),
\end{split}
\end{equation}
such 
that for an analytic function $\mathcal{F}$ applied to $A$, we obtain
\begin{equation}\label{eqn:DiagFLaplace2D}
\mathcal{F}(A)  = (F_1^*\otimes F_2^*) 	\mathcal{F}(\varLambda)	(F_1\otimes F_2 ).
\end{equation}

Now assume that the diagonal matrix $\mathcal{F}(\varLambda)$ can be expressed approximately 
as a short sum of Kronecker rank-$1$ diagonal matrices, i.e.,
\begin{equation*}
	\mathcal{F}(\varLambda)  = \sum_{k=1}^R \diag \big(\mathbf{u}_1^{(k)} \otimes \mathbf{u}_2^{(k)}\big),
\end{equation*}
with vectors $\mathbf{u}_\ell^{(k)}\in\R^{n_\ell}$ and $R \ll \min(n_1,n_2)$, and, moreover,
let
$\mathbf{x}\in \R^N$ be a vector given in a low rank Kronecker form, i.\,e.
\begin{equation*}
	\mathbf{x} = \sum_{j=1}^S \mathbf{x}_1^{(j)} \otimes \mathbf{x}_2^{(j)},
\end{equation*}
with vectors $\mathbf{x}_\ell^{(j)}\in\R^{n_\ell}$ and $S \ll \min(n_1,n_2)$. Then, 
we can compute a matrix-vector product
\begin{equation}\label{eqn:Lap2DtimesX}
\begin{split}
	\mathcal{F}(A) \mathbf{x} &= (F_1^*\otimes F_2^*) 
	\bigg( \sum_{k=1}^R \diag \big(\mathbf{u}_1^{(k)} \otimes \mathbf{u}_2^{(k)}\big) \bigg)
	(F_1\otimes F_2 )
	\bigg( \sum_{j=1}^S \mathbf{x}_1^{(j)} \otimes \mathbf{x}_2^{(j)} \bigg)\\
	&= \sum_{k=1}^R \sum_{j=1}^S F_1^* \big( \mathbf{u}_1^{(k)} \odot F_1 \mathbf{x}_1^{(j)} \big) \otimes
	F_2^*\big( \mathbf{u}_2^{(k)} \odot F_2 \mathbf{x}_2^{(j)} \big),
\end{split}
\end{equation}
where $\odot$ denotes the Hadamard (componentwise) product.
Using the sin-FFT, expression \eqref{eqn:Lap2DtimesX} can be calculated in factored 
form in $\mathcal{O}(RSn\log n)$ flops, where $n = \max(n_1,n_2)$.

Likewise, in the case $d=3$, by completely analogous reasoning, equation \eqref{eqn:Lap2DtimesX} becomes
\begin{equation}\label{eqn:Lap3DtimesX}
	\mathcal{F}(A) \mathbf{x} = 
	\sum_{k=1}^R \sum_{j=1}^S F_1^* \big( \mathbf{u}_1^{(k)} \odot F_1 \mathbf{x}_1^{(j)} \big) \otimes
	F_2^*\big( \mathbf{u}_2^{(k)} \odot F_2 \mathbf{x}_2^{(j)} \big)
	\otimes F_3^*\big( \mathbf{u}_3^{(k)} \odot F_3 \mathbf{x}_3^{(j)} \big),
\end{equation}
and similar in the case of $d>3$, which ensure the evaluation cost $\mathcal{O}(d RSn\log n)$.

In our application the matrix function $F(A)$  has the form 
$F(A)=\big(\beta A^{-\alpha} + \tfrac{\gamma}{\beta}A^\alpha\big)$ and its inverse, $(F(A))^{-1}$.
In what follows, we discuss the tensor method for solving the linear system of equations (\ref{eqn:forControl})
in the low-rank formats. To that end, we use the PCG scheme on the ``low-rank'' manifold 
sketched  in the next paragraph.
This algorithm allows the low cost calculations with rank truncations applied to
both 2D and 3D cases. 
In the 3D case, we use the rank truncation in the canonical format 
by application of the RHOSVD decomposition (see Section \ref{sec:low-rank_tensor}).

\subsection{The low-rank PCG scheme}\label{ssec:Rank_PCG}

For operators $\mathtt{func}$ and $\mathtt{precond}$ given in a low-rank format, 
such as \eqref{eqn:Lap2DtimesX} (for $d=2$) or \eqref{eqn:Lap3DtimesX} (for $d=3$), Krylov subspace methods
can be applied very efficiently, since they only require matrix-vector products. The formulation
of the PCG method in Algorithm~\ref{alg:pcg} is independent of $d$, 
as long as an appropriate rank truncation procedure $\mathtt{trunc}$ is chosen.

\begin{algorithm}[H]
	\caption{Preconditioned CG method in low-rank format} \label{alg:pcg}
	\begin{algorithmic}[1]
		\Require{Rank truncation procedure $\mathtt{trunc}$, rank tolerance parameter $\varepsilon$,
		linear function in low rank format $\mathtt{fun}$, preconditioner in low rank format
		$\mathtt{precond}$, right-hand side tensor $\mathbf{B}$, initial guess $\mathbf{X}^{(0)}$,
		maximal iteration number $k_{\max}$}
			\State $\mathbf{R}^{(0)} \leftarrow \mathbf{B} - \mathtt{fun}(\mathbf{X}^{(0)})$
			\State $\mathbf{Z}^{(0)} \leftarrow \mathtt{precond}(\mathbf{R}^{(0)})$
			\State $\mathbf{Z}^{(0)} \leftarrow \mathtt{trunc}(\mathbf{Z}^{(0)},\varepsilon)$
			\State $\mathbf{P}^{(0)} \leftarrow (\mathbf{Z}^{(0)})$
			\State $k \leftarrow 0$
			\Repeat
				\State $\mathbf{S}^{(k)} \leftarrow
				\mathtt{fun}(\mathbf{P}^{(k)})$
				\State {\color{red}$\mathbf{S}^{(k)} \leftarrow
				\mathtt{trunc}(\mathbf{S}^{(k)},\varepsilon)$}
				\State $\alpha_k \leftarrow
				\frac{\dotprod{\mathbf{R}^{(k)}}{\mathbf{Z}^{(k)}}}
				{\dotprod{\mathbf{P}^{(k)}}{\mathbf{S}^{(k)}}}$
				\State $\mathbf{X}^{(k+1)} \leftarrow
				\mathbf{X}^{(k)} + \alpha_k \mathbf{P}^{(k)}$
				\State {\color{red}$\mathbf{X}^{(k+1)} \leftarrow
				\mathtt{trunc}(\mathbf{X}^{(k+1)},\varepsilon)$}
				\State $\mathbf{R}^{(k+1)} \leftarrow
				\mathbf{R}^{(k)} - \alpha_k \mathbf{S}^{(k)}$
				\State {\color{red}$\mathbf{R}^{(k+1)} \leftarrow
				\mathtt{trunc}(\mathbf{R}^{(k+1)},\varepsilon)$}
				\If {$\mathbf{R}^{(k+1)}$ is sufficiently small}
					\State \Return $\mathbf{X}^{(k+1)}$
					\State \textbf{break}
				\EndIf
				\State $\mathbf{Z}^{(k+1)} \leftarrow
				\mathtt{precond}(\mathbf{R}^{(k+1)})$
				\State {\color{red}$\mathbf{Z}^{(k+1)} \leftarrow
				\mathtt{trunc}(\mathbf{Z}^{(k+1)},\varepsilon)$}
				\State $\beta_k \leftarrow
				\frac{\dotprod{\mathbf{R}^{(k+1)}}{\mathbf{Z}^{(k+1)}}}
				{\dotprod{\mathbf{Z}^{(k)}}{\mathbf{R}^{(k)}}}$
				\State $\mathbf{P}^{(k+1)} \leftarrow
				\mathbf{Z}^{(k+1)} + \beta_k \mathbf{P}^{(k)}$
				\State {\color{red} $\mathbf{P}^{(k+1)} \leftarrow
				\mathtt{trunc}(\mathbf{P}^{(k+1)},\varepsilon)$}
				\State $k \leftarrow k+1$
			\Until{$k = k_{\max}$}
			\Ensure{Solution $\mathbf{X}$ of $\mathtt{fun}(\mathbf{X})=\mathbf{B}$}
	\end{algorithmic}
\end{algorithm}
As the rank truncation procedure,
in our implementation we apply the reduced SVD algorithm in 2D case and the RHOSVD based 
canonical-to-Tucker-to-canonical algorithm (see \cite{KhKh3:08}) 
as described in Section \ref{sec:low-rank_tensor}.

\section{Low-rank tensor approximation for analytic functions of fractional $d$-dimensional Laplacian}
\label{sec:tensor_approx_theory}

\subsection{Rank-structured decomposition for functions of the core diagonal matrix $\Lambda$}
\label{ssec:Multivariate_cores}

In this section we analyze the rank-structured tensor decompositions of various matrix (tensor) 
valued functions on the discrete Laplacian arising in the solution of problems 
(\ref{eqn:Lagrange_cont}) and (\ref{eqn:State}) including different combinations of fractional 
Laplacian  in $\mathbb{R}^d$. These decompositions provide the short term Kronecker product representations 
for functions of multidimensional Laplacian.

We consider the matrices $A_1, A_2$ and $A_3$ defined as the matrix valued functions of 
the discrete Laplacian ${A}_h$ by the equations
\begin{equation} \label{eqn:L_m_alpha}
{A}_1= {A}_h^{-\alpha} ,
\end{equation}
\begin{equation} \label{eqn:L_alp_m_alp}
{A}_2= {A}_h^{-\alpha} + {A}_h^{\alpha},
\end{equation}
and 
\begin{equation} \label{eqn:(L_alp_m_alp)_m1}
 {A}_3= \big({A}_h^{-\alpha} + {A}_h^{\alpha}\big)^{-1}=A_2^{-1},
\end{equation}
respectively. It is worth to notice that the matrix $A_3$ defines the solution operator
in equation (\ref{eqn:Lagrange_cont}), which allows to calculate the optimal control
in terms of the design function ${\bf y}_\Omega$ on the right-hand side of (\ref{eqn:Lagrange_cont})
 by solving the equation 
\begin{equation} \label{eqn:OptCont}
 {A}_2 {\bf u}^\ast=  {\bf y}_\Omega. 
\end{equation}
Finally, the state variable is calculated by 
\begin{equation} \label{eqn:L_alp_state}
{\bf y}= \beta {A}_h^{-\alpha} {\bf u}^\ast = \beta A_1 {\bf u}^\ast.
\end{equation}

In the presented numerics we consider the rank bounds of the Tucker/canonical (or SVD in the 2D case) 
decompositions
for the corresponding multi-indexed core tensors/matrices $\mathcal{F}(\varLambda_p)$, $p=1,2,3$,
further denoted by $G_1,G_2,G_3\in \mathbb{R}^{n\times n}$ in the 2D case, and by 
${\bf G}_1,{\bf G}_2,{\bf G}_3\in \mathbb{R}^{n\times n\times n}$ in the 3D case
and representing the matrix valued functions $A_1,A_2,A_3$ of ${A}_h$  in the Fourier basis, see 
(\ref{eqn:DiagFLaplace3D}) and (\ref{eqn:DiagFLaplace2D}).
This factorization is well suited for the rank-structured algebraic operations since 
the $d$-dimensional Fourier transform matrix has the Kronecker rank equals to one. 
For example, for $d=3$ we have
\[
 {\cal F}=F_1 \otimes F_2 \otimes F_3.
\]
Let $\{\lambda_i\}_{i=1}^n$ be the set of eigenvalues for the 1D finite difference Laplacian
in $H^{1}_0(0,1)$ discretized on the uniform grid with the mesh size $h=1/(n+1)$, see (\ref{eqn:eig_lap}).
The elements of the core diagonal matrix $\varLambda$ in (\ref{eqn:DiagLaplace}) 
can be reshaped to a three-tensor 
$$
{\bf G}=[g(i_1,i_2,i_3)]\in \mathbb{R}^{n_1 \times n_2 \times n_3}, 
\quad i_\ell\in \{1,\ldots,n_\ell\},
$$
where
\[
 g(i_1,i_2,i_3)=\lambda_{i_1} + \lambda_{i_2} + \lambda_{i_3},
\]
implying that the three-tensor ${\bf G}$ has the exact rank-$3$ canonical decomposition. In the case $d=2$, we have 
similar two-term representation for the matrices $G=[g(i_1,i_2)]\in \mathbb{R}^{n_1 \times n_2}$, 
$g(i_1,i_2)=\lambda_{i_1} + \lambda_{i_2} $.

 In the 3D case we consider the Tucker and canonical decompositions of the $n\times n \times n$ core 
 tensors,  corresponding to the matrices $A_1,A_2,A_3$, 
\begin{equation}\label{eq:coreTens} 
 {\bf G}_p=[g_p(i,j,k)], \quad p=1,2,3 
\end{equation}
with entries defined by
 \begin{equation}  \label{eq:3a}
 g_1 (i,j,k) = \frac{1}{(\lambda_i +\lambda_j +\lambda_k)^\alpha},
 \end{equation} 
\begin{equation}\label{eq:3b}
 g_2 (i,j,k) = \frac{1}{(\lambda_i +\lambda_j +\lambda_k)^\alpha} + 
 (\lambda_i +\lambda_j+\lambda_k)^\alpha,
\end{equation}
 \begin{equation}\label{eq:3c}
 g_3 (i,j,k) = \left((\lambda_i +\lambda_j+\lambda_k)^{-\alpha} + 
 (\lambda_i +\lambda_j+\lambda_k)^\alpha \right)^{-1}.
\end{equation}

In the 2D case we analyze the singular value decomposition of the $n\times n $ core matrices 
\begin{equation}\label{eq:coreMatr}
 G_p=[g_p(i,j)], \quad   p=1,\, 2,\, 3,  
\end{equation}
with entries defined by
\begin{equation} \label{eq:2a}
 g_1 (i,j) = \frac{1}{(\lambda_i +\lambda_j)^\alpha},
 \end{equation} 
\begin{equation}\label{eq:2b}
 g_2 (i,j) = \frac{1}{(\lambda_i +\lambda_j)^\alpha} + (\lambda_i +\lambda_j)^\alpha,
\end{equation}
 \begin{equation}\label{eq:2c}
 g_3 (i,j) = \left((\lambda_i +\lambda_j)^{-\alpha} + (\lambda_i +\lambda_j)^\alpha \right)^{-1}.
\end{equation}

The error estimate for the rank decomposition of the matrices 
$G_p$ and the respective 3D tensors ${\bf G}_p$ can be derived based on the sinc-approximation 
theory as discussed in \S\ref{ssec:Theory_multivariate_cores}. 
We consider the class of matrix valued functions of the discrete Laplacian, 
$A_1,\ldots,A_3$, given 
by (\ref{eqn:L_m_alpha}) -- (\ref{eqn:(L_alp_m_alp)_m1}).
In view of the FFT diagonalization, the tensor approximation problem is reduced to the analysis 
of the corresponding function related tensors ${\bf G}_1,\ldots,{\bf G}_3$ specified by 
multivariate functions of the discrete argument, $g_1,\ldots,g_3$, given by (\ref{eq:3a}) -- (\ref{eq:3c}).

 It is worse to note that the sinc quadrature approximation theory based on the Laplace transform applies only to 
the operators (matrices) with the negative fractional power like ${{A}_h}^{-\alpha}$, $\alpha >0$.
To prove the existence of rank decomposition for the positive power of ${A}_h$, we notice that for 
$1>\alpha >0$
\[
 A^{\alpha} = A \cdot A^{\alpha-1}, \quad \alpha-1 <0,
\]
where the second factor on the right-hand side can be approximated with the low Kronecker 
rank due to the general theory of sinc approximation,
while the initial stiffness matrix $A$ has the Kronecker rank equals to $d$. 
Similar argument applies to the case $\alpha >1$. This proves the following Lemma.
\begin{lemma}\label{lem:Rank_positive_power}
For given $\alpha >0$, let the matrix $A^{-\beta}$, $\beta>0$, can be represented in the Kronecker rank-$R$ 
form. Then the matrix $A^{\alpha}$ can be represented with the Kronecker rank at most $d\, R$.
\end{lemma}

\begin{remark}\label{rem:direct_fract_Lap}
Our numerical scheme can be also applied to the class of ``directionally fractional'' Laplacian-type operators 
${\cal A}_{(\alpha)}$, obtained by the essential simplification of the fractional 
elliptic operator ${\cal A}^\alpha$ in (\ref{eqn:frac_setting_Lap}) considered 
in our paper,  as follows
\begin{equation}  \label{eq:Laplace_1Dfract}
{\cal A}_{(\alpha)} := 
\sum\limits_{\ell=1}^d 
\left(-\Delta_\ell \right)^\alpha, \quad \alpha>0,
\end{equation} 
where $\Delta_\ell$ is the 1D Laplacian in variable $x_\ell$.
In the case $d=3$ the core tensors $g_p$ in (\ref{eq:3a}) -- (\ref{eq:3c}) representing the operator 
${\cal A}_{(\alpha)}$ in the Fourier basis are simplified to 
 \begin{equation}  \label{eq:3a1d}
 \widetilde{g}_1 (i,j,k) = \frac{1}{\lambda_i^\alpha +\lambda_j^\alpha +\lambda_k^\alpha},
 \end{equation} 
\begin{equation}\label{eq:3b1d}
 \widetilde{g}_2 (i,j,k) = \frac{1}{\lambda_i^\alpha +\lambda_j^\alpha +\lambda_k^\alpha} +
 \lambda_i^\alpha +\lambda_j^\alpha +\lambda_k^\alpha,
  \end{equation}
 \begin{equation}\label{eq:3c1d}
 \widetilde{g}_3 (i,j,k) = \left(\frac{1}{\lambda_i^\alpha +\lambda_j^\alpha +\lambda_k^\alpha} +
 \lambda_i^\alpha +\lambda_j^\alpha +\lambda_k^\alpha \right)^{-1},
\end{equation}
respectively. 
The low rank tensor decomposition of these discrete functions 
is practically identical to the case of Laplacian operator in (\ref{eq:3a}) -- (\ref{eq:3c})
corresponding to the fractional power $\alpha=1$. 
The similar model has been considered in \cite{DPSS:2014} for $d=2,3$, 
where instead of $(-\Delta_\ell)^{\alpha}$ the so-called Caputo type integral representation has been 
adapted. Instead of the Fourier based diagonalization, the QTT tensor approximation 
\cite{KhQuant:09,Osel-TT-LOG:09} of the arising $n\times n$ univariate Toeplitz matrices has been applied.
\end{remark}

The application of the tensor techniques to the case of variable coefficients in 
(\ref{eq:Laplace_1Dfract}) in 3D case requires the solution of the few 1D spectral problems 
for symmetric three-diagonal $n\times n$ matrices, which on the order of magnitude faster 
compared with the optimal cost $O(n^4)$ for solving the 3D control problem in the full matrix format.

\begin{remark}\label{rem:generalized_fract_Lap}
In the general case of ``directionally fractional'' operators with variable 
coefficients $a_\ell(x_\ell)>0$ in (\ref{eq:Laplace_1Dfract})
the orthogonal matrices $F_\ell$, $\ell=1,\ldots,d$, 
of the univariate Fourier transforms in (\ref{eqn:DiagLaplace})
should be substituted by the orthogonal matrices of the eigenvalue decomposition for the 
discretized elliptic operators $A_{h,\ell}$ (1D stiffness matrices) corresponding to the substitution 
$-\Delta_\ell  \mapsto A_\ell= -\nabla_\ell^T \cdot a_\ell(x_\ell)\nabla_\ell$. The eigenvalues in 
(\ref{eq:3a1d}) -- (\ref{eq:3c1d}) are obtained from the solution of the discrete eigenvalue problem 
$A_{h,\ell} {\bf u}_i = \lambda_{i,\ell} {\bf u}_i$, $\ell=1,\ldots,d$.
\end{remark}

\subsection{Rank estimates for matrix valued functions}
\label{ssec:Theory_multivariate_cores}

In this section, we sketch the proof of the existence of the low 
rank canonical/Tucker decomposition of the core tensors ${\bf G}_p$, $p=1,\ldots,4$.
Following \cite{Stenger}, we
define the {Hardy space} $H^1(D_{\delta})$ of functions which 
are analytic in the strip
\[
D_{\delta}:=\{z\in\mathbb{C}:|\Im \,z|\leq\delta\},\quad 0<\delta<\frac{\pi}{2},
\]
and satisfy
\[ 
N(f,D_{\delta}):=\int_{\mathbb{R}}\left(\left\vert f(x+i\delta)\right\vert
  +\left\vert f(x-i\delta)\right\vert \right)  dx<\infty. 
\]
Recall that for $f\in H^1(D_{\delta})$, the integral
\begin{equation*} 
I(f)=\int_{\Omega}f(x)d x \quad\left(  \Omega=\mathbb{R}\;\;\text{or}%
  \;\Omega=\mathbb{R}_{+} \right)   \label{(I(f)}%
\end{equation*}
can be approximated by the {Sinc-quadrature} (trapezoidal rule) 
\[
T(f,\mathfrak{h}):=\mathfrak{h}\sum_{k=-\infty}^{\infty}f(t_k),\;\; t_k=k\mathfrak{h};
\quad \quad
|I(f) - T(f,\mathfrak{h})|=O(e^{-\pi \delta/\mathfrak{h}}),
\quad \mathfrak{h}\to 0,
\]
that converges exponentially fast in $\mathfrak{h}\to 0$.
Similar estimates hold for computable truncated sums, see \cite{Stenger}
\begin{equation}\label{eqn:Sinc_quadr_trunc}
T_{M}(f,\mathfrak{h}) :=  \mathfrak{h}\sum_{k=-M}^{M}f(k\mathfrak{h}).
\end{equation}
Indeed, if $f\in {H}^{1}(D_{\delta})$ and 
$ | f(x)|\leq C\exp(-b|x|)$ for all $x \in\mathbb{R}$ $\;b,C>0, $ then 
\begin{equation} 
|  I(f)-T_{M}(f,\mathfrak{h})| \leq C\left[  \frac{e^{-2\pi 
\delta/\mathfrak{h}}}{1-e^{-2\pi\delta/\mathfrak{h}}}N(f,D_{\delta})+\frac 
{1}{b}e^{-b\mathfrak{h}M}\right].  \label{ff6}%
\end{equation}

In our context, the Sinc-quadrature approximation applies to multivariate functions
$F:\mathbb{R}^d \to \mathbb{R}$ of a sum of single variables, say,
$F(x_1,\ldots,x_d)=f(\rho)$ with $\rho=\sum_{\ell=1}^{d} f_\ell(x_{\ell}) >0$,
where $f_\ell: \mathbb{R} \to \mathbb{R}_+$, 
by using the integral representation of the analytic univariate function 
$f:\mathbb{R}_+\to \mathbb{R}$, 
\[ 
f(\rho)=\int_{\Omega} G(t) e^{-\rho E(t)}\mathrm{d}t\approx \sum_{k=1}^R c_k e^{-\rho E(t_k)} 
= \sum_{k=1}^R c_k  \prod_{\ell=1}^{d} e^{- f_\ell (x_\ell) E(t_k)}, \quad 
 \Omega\in \{ \mathbb{R},\mathbb{R_+}\}.
\] 
 In the cases
(\ref{eqn:L_m_alpha}) -- (\ref{eqn:(L_alp_m_alp)_m1}) and (\ref{eqn:L_alp_state})
the related functions $f(\rho)$ take the particular form
\[
f(\rho)=\rho^{\alpha}, \quad f(\rho)=\rho^{-\alpha}, \quad 
f(\rho)=(\rho^{\alpha} +\rho^{-\alpha})^{-1}. 
\]
The univariate function $f$ may have point singularities or cusps at $\rho=0$, 
say, $f(\rho)=\rho^{\pm \alpha}$.
Applying the Sinc-quadrature (\ref{eqn:Sinc_quadr_trunc}) to the Laplace-type  transform leads to 
{rank-$R$ separable approximation} of the function $F$,
\begin{equation}\label{eqn:Func_of_rho}
F(x)=f(\rho)= f(f_1(x_{1})+\ldots+  f_d(x_{d}))\approx
\sum_{k=1}^R\omega_{k}G(t_{k})  e^{-\rho E(t_{k})}%
=
\sum\limits_{k=1}^R  c_{k}\prod_{\ell=1}^{d} e^{-f_\ell(x_{\ell})E(t_{k})}, 
\end{equation}
with $c_{k}=\omega_{k}G(t_{k})$ and $R=2M+1$. 

Notice that the generating function $f$ can be defined 
on the discrete argument, i.e., on the multivariate index 
${\bf i}=(i_1,\ldots,i_d)$, $i_\ell=1,\ldots,n$,
such that each univariate function $f_\ell$ in (\ref{eqn:Func_of_rho})
is defined on the index set $\{i_\ell\}$, $\ell=1,\ldots ,d$.
In our particular applications to functions of the discrete fractional 
Laplacian ${A}_h^{-\alpha}$ we have
\begin{equation}\label{eqn:eig_lap1}
f_\ell(i_\ell)=\lambda_{i_\ell} = -\frac{4}{h_\ell^2}\sin^2 \bigg( \frac{\pi i_\ell h_\ell}{2} \bigg),
\end{equation}
where $\lambda_{i_\ell}$ denote the eigenvalues of the univariate discrete Laplacian with the Dirichlet 
boundary conditions, see (\ref{eqn:eig_lap}) and (\ref{eq:3a}) -- (\ref{eq:3c}).

In this case the integral representation of the function $f(\rho)=\rho^{-\alpha}$, 
$\rho = \sum_{\ell=1}^{d} f_\ell(x_{\ell}) >0$, with $f_\ell$ given by (\ref{eqn:eig_lap1}) 
$\alpha >0$, takes a form
\begin{equation} \label{eqn:Laplace_transf_rhom1}
 {\rho}^{-\alpha}=\frac{1}{\Gamma(\alpha)} \int_0^\infty t^{\alpha -1} e^{- {\rho} t}\, dt,
 \quad \rho\in [1,B], \; B > 1.
\end{equation}


Several efficient sinc-approximation schemes for classes of multivariate functions and operators 
have been developed, see
\cite{GHK:05,HaKhtens:04I,Khor1:06,KhorBook:17}. In the case (\ref{eqn:Laplace_transf_rhom1})
the simple modification of Lemma 2.51 in \cite{KhorBook:17} can be applied, see also \cite{HaKhtens:04I}.
To that end, the substitution $t= \phi(u):=\log(1+e^u)$, that maps $\phi:\mathbb{R} \to\mathbb{R}_+ $,
leads to the integral
\[
 {\rho}^{-\alpha}=\frac{1}{\Gamma(\alpha)} \int_\mathbb{R} f_1(u) du , \quad 
 f_1(u)= \frac{[\log(1+e^u)]^{\alpha-1} e^{-\rho \log(1+e^{-u})}}{1+e^{-u}},
\]
which can be approximated by the sinc-quadrature with the choice 
$\mathfrak{h}=\frac{\gamma \pi}{\sqrt{M}}$ with some $ 0 < \gamma < 1$.
This argument justifies the accurate low-rank representation of functions 
in (\ref{eq:2a}) and (\ref{eq:3a}) representing the fractional Laplacian inverse.

\subsection{On numerical validation and some generalizations}
\label{ssec:Numerics_generalization}

The numerical results presented in Section \ref{sec:numerics_tensor} clearly 
illustrate the high accuracy of  the low-rank approximations to  
various matrix valued functions of the fractional Laplacian defined in (\ref{eq:2a}) -- (\ref{eq:2c})
  and in (\ref{eq:3a}) -- (\ref{eq:3c}). In our numerical tests the rank decomposition 
  of the 2D and 3D tensors under consideration was performed by the multigrid Tucker-to-canonical 
  scheme as described in Section \ref{sec:low-rank_tensor}.

 \begin{remark}\label{rem:GenerLapl}
 The presented approach is applicable with minor modifications to the case of 
 more general elliptic operators of the form (for $d=3$)
 \begin{equation} \label{eqn:Laplace_modif}
  {\cal A}= {F}(-\Delta_1)\otimes I_2 \otimes I_3 + I_1\otimes {F}(-\Delta_1)\otimes I_3
  + I_1\otimes I_2 \otimes {F}(-\Delta_1),
 \end{equation}
 where $F$ is the rather general analytic function of the univariate Laplacian.  
 In this case the fractional operators ${\cal A}^{\alpha}$ and ${\cal A}^{-\alpha}$
 can be introduced by the similar way, where the values $\{\lambda_i\}$ in the representation of the
 core tensors ${\bf G}_p$, $p=1,2,3$, should be substituted by $\{F(\lambda_i)\}$. 
 Given the symmetric positive definite matrix $X\in \mathbb{R}^{n\times n}$, 
 then the particular  choice $F(X)=X^{\pm \alpha}$, $\alpha>0$, suites well for the presented approach.
  In general, the cost of matrix-vector multiplication for matrices ${F}(X)$ is estimated by $O(n^2)$.
 However, this cost can be reduced to the logarithmic scale by using the QTT tensor format 
 \cite{KhQuant:09,Osel-TT-LOG:09,KhorBook:17}.     
\end{remark}
 
 We notice that the rank structured approximation of the control problem with 
 the Laplacian type operator in the form (\ref{eqn:Laplace_modif}) with the particular choice 
 $F(-\Delta_1)=(-\Delta_1)^{\alpha}$ was considered in \cite{DPSS:2014}.
 Since the operator in (\ref{eqn:Laplace_modif}) only includes the univariate 
 fractional Laplacain, this case can be treated as for the standard 
 Laplacian type control problems as pointed out in Remark \ref{rem:GenerLapl}.

\section{Basics of tensor decomposition for discretized multivariate functions and operators}
\label{sec:low-rank_tensor}

Here we recall the rank-structured tensor decompositions and  tensor numerical techniques used in this paper.
The basic rank-structured representations of the multidimensional tensors are the 
canonical \cite{Hitch:1927} and Tucker \cite{Tuck:1966} tensor formats. 
They have been since long used in computer science for quantitative analysis of data 
in chemometrics, psychometrics  and signal processing \cite{DMV-SIAM2:00}.
In the last decades an extensive research have been focused on different aspects 
of multilinear algebra, tensor numerical calculus and related issues, 
see \cite{DMV-SIAM2:00,Cichocki:2002,Hack_Book12,KhorBook:17,Khor2Book2:18,MaRaSch:19} for an overview. 
 
A tensor of order $d$ in a full format,  is defined as a 
multidimensional array over a $d$-tuple index set:
\begin{equation} \label{Tensor_def}
{\bf T}=[t_{i_1,\ldots,i_d}] 
\in \mathbb{R}^{n_1 \times \ldots \times n_d}\,
\end{equation}
with indices $i_\ell\in I_\ell:=\{1,\ldots,n_\ell\}, \; \ell=1,\ldots d $. 
For a tensor ${\bf T}$ given entry-wise in a form (\ref{Tensor_def}), 
(assuming $n_\ell =n$) both the  required storage  
and complexity of operations scale   exponentially
in the dimension size $d$, as $O(n^d)$, \cite{bellman-dyn-program-1957}. 

To avoid the exponential scaling in the dimension, the 
rank-structured separable approximations  of the multidimensional 
tensors can be used.  
A tensor in the $R$-term canonical format is defined by a sum of rank-$1$ 
tensors 
 \begin{equation}\label{eqn:CP_form}
   {\bf T} = {\sum}_{k =1}^{R} \xi_k
   {\bf u}_k^{(1)}  \otimes \ldots \otimes {\bf u}_k^{(d)},  \quad  \xi_k \in \mathbb{R},
\end{equation}
where ${\bf u}_k^{(\ell)}\in \mathbb{R}^{n_\ell}$ are normalized vectors, 
and $R$ is the canonical rank. The storage cost of this
parametrization is bounded by  $d R n$.
An alternative (contracted product) notation for a canonical tensor can be 
used (cf. the Tucker tensor format, (\ref{eqn:Tucker_form}))
\begin{equation}
 \label{Can_Contr}
 {\bf T} = {\bf C} \times_1 U^{(1)} \times_2 U^{(2)} \times_3 \cdots \times_d U^{(d)},
 \end{equation}
where ${\bf C}=diag\{\xi_1,...,\xi_R\} \in\mathbb{R}^{R^{\otimes d}}$ is a super-diagonal tensor,  
$ U^{(\ell)}=[{\bf u}_1^{(\ell)}\ldots {\bf u}_R^{(\ell)} ]\in\mathbb{R}^{n_\ell \times R}$ are
the side matrices, and $\times_\ell$ denotes the contracted product in mode $\ell$.
  For $d\geq 3$, computation of the  canonical rank representation (\ref{eqn:CP_form}) 
for a tensor ${\bf T}$ in a form (\ref{Tensor_def}) is, in general, an $N$-$P$ hard problem.

The orthogonal Tucker tensor format is suitable for stable numerical decompositions with a fixed
truncation threshold.
We say that the tensor ${\bf T} $ is represented in the rank-$\bf r$ orthogonal Tucker format 
with the rank parameter ${\bf r}=(r_1,\ldots,r_d)$ if 
\begin{equation}
\label{eqn:Tucker_form}
  {\bf T}  =\sum\limits_{\nu_1 =1}^{r_1}\ldots
\sum\limits^{r_d}_{{\nu_d}=1}  \beta_{\nu_1, \ldots ,\nu_d}
\,  {\bf v}^{(1)}_{\nu_1} \otimes  
{\bf v}^{(2)}_{\nu_2} \ldots \otimes {\bf v}^{(d)}_{\nu_d}, 
\end{equation}
where $\{{\bf v}^{(\ell)}_{\nu_\ell}\}_{\nu_\ell=1}^{r_\ell}\in \mathbb{R}^{n_\ell}$, $\ell=1,\ldots,d$
represents a set of orthonormal vectors 
and $\boldsymbol{\beta}=[\beta_{\nu_1,\ldots,\nu_d}] \in \mathbb{R}^{r_1\times \cdots \times r_d}$ is 
the Tucker core tensor. 
The storage cost for the Tucker tensor format is bounded by $d r n +r^d$, 
with $r=|{\bf r}|:=\max_\ell r_\ell$.
Using the orthogonal side matrices $V^{(\ell)}=[{\bf v}^{(\ell)}_1 \ldots {\bf v}^{(\ell)}_{r_\ell}]$, 
the Tucker tensor decomposition can be presented by using
contracted products, 
\begin{equation}\label{Tuckt} 
  { {\bf T}_{({\bf r})}=\boldsymbol{\beta} \times_1 V^{(1)}\times_2 V^{(2)} \times_3 \ldots 
   \times_d V^{(d)}.} 
\end{equation}  
 The problem of the best Tucker tensor approximation to a given tensor 
  ${\bf T}_0\in \mathbb{R}^{n_1 \times \ldots \times n_d}$ is 
the following minimization problem, 
\begin{equation} \label{Tuck Cost f1}
 f({\bf T}):= \| {\bf T} - {\bf T}_0 \|^2\to min
\quad\mbox{over}\;\; {\bf T} \in \mbox{\boldmath{$\mathcal{T}$}}_{{\bf r},{\bf n}},
\end{equation}
in the class of rank-${\bf r}$ Tucker tensors $\mbox{\boldmath{$\mathcal{T}$}}_{{\bf r},{\bf n}}$, 
that is equivalent to the maximization problem \cite{DMV-SIAM2:00}
\begin{equation}\label{GRQ norm r}
g(V^{(1)},...,V^{(d)}) :=
\left\|{\bf T}_0 \times_1 {V^{(1)}}^T\times ...
  \times_d {V^{(d)}}^T \right\|^2 \to max
\end{equation}
over the product set of ${\cal M}_\ell$ of orthogonal matrices
$V^{(\ell)}\in \mathbb{R}^{n_\ell \times r_\ell} $, $\ell=1,2,\ldots, d$.
For given maximizers $V^{(\ell)} $, the core $\boldsymbol{\beta}$ 
that minimizes (\ref{Tuck Cost f1}) is computed as
\begin{equation}\label{coret}
\boldsymbol{\beta}= {\bf T}_0\times_1 {V^{(1)}}^T\times_2 \ldots
  \times_d {V^{(d)}}^T \in \mathbb{R}^{r_1\times ...\times r_d},
   \end{equation}  
yielding contracted product representation (\ref{Tuckt}). 
  The Tucker tensor format provides a stable decomposition algorithm \cite{DMV-SIAM2:00}, which requires 
the tensor in a full size format $O(n^d)$. The main step of this algorithm, the higher order singular 
value decomposition (HOSVD), 
features the complexity of the order of $O(n^{d+1})$. This restricts the applicability of 
the algorithm to moderate size  tensors.

For the wide class of function related tensors the Tucker tensor decomposition  
exhibits exceptional approximation properties 
(logarithmically low ranks) \cite{Khor1:06,KhKh:06}. 
This enables the multigrid Tucker decomposition for full size tensors
and the  reduced higher order SVD (RHOSVD) for the Tucker decomposition of tensors given 
in the canonical format (but with possibly large rank). 
 The RHOSVD provides linear complexity in the univariate 
grid-size $n$, $O(dnR^2)$, independent on dimension $d$, since it does not require the construction of the 
full size tensor  \cite{KhorBook:17,Khor2Book2:18}.

At the intermediate steps of the algorithms we use the mixed Tucker-canonical transform \cite{KhKh3:08}.
Given the rank parameters ${\bf r},R$, the explicit representation of tensor ${\bf T}$ 
in a mixed Tucker-canonical tensor format with the core $\boldsymbol{\beta}$ represented in 
the canonical format is given by
\begin{equation}\label{two_level}
 {\bf T} = \left( \sum_{\nu =1}^R \xi_\nu
   {\bf u}^{(1)}_\nu  \otimes \ldots \otimes {\bf u}^{(d)}_\nu
 \right)  \times_1 V^{(1)}\times_2 V^{(2)} \ldots \times_d V^{(d)}.
\end{equation}
The corresponding side-matrices for the resulting canonical tensor
are given by $U^{(\ell)}=[V^{(\ell)}{\bf u}^{(\ell)}_1 \ldots V^{(\ell)}{\bf u}^{(\ell)}_{R}]$,
with the scaling coefficients  $\xi_{\nu}$ ($\nu=1,\ldots ,R$), where $R \leq r^2$ with $r=\max r_\ell$ 
holds for $d=3$, see \cite{Khor2Book2:18} for further details. The mixed
tensor format is efficient for approximation of function related tensors, where
the exponentially fast  decay in $\bar{r}=\min r_\ell$ of the Tucker approximation error   
leads to small sizes of the Tucker core \cite{Khor1:06,KhKh:06}.

\begin{figure}[tbh]
\centering
\includegraphics[width=7.0cm]{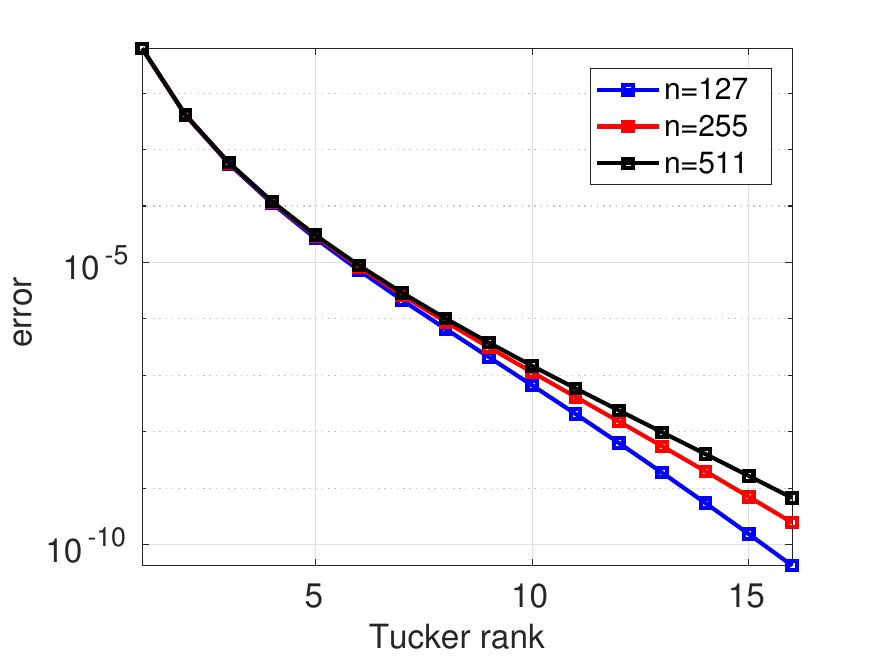}\quad
 \includegraphics[width=7.0cm]{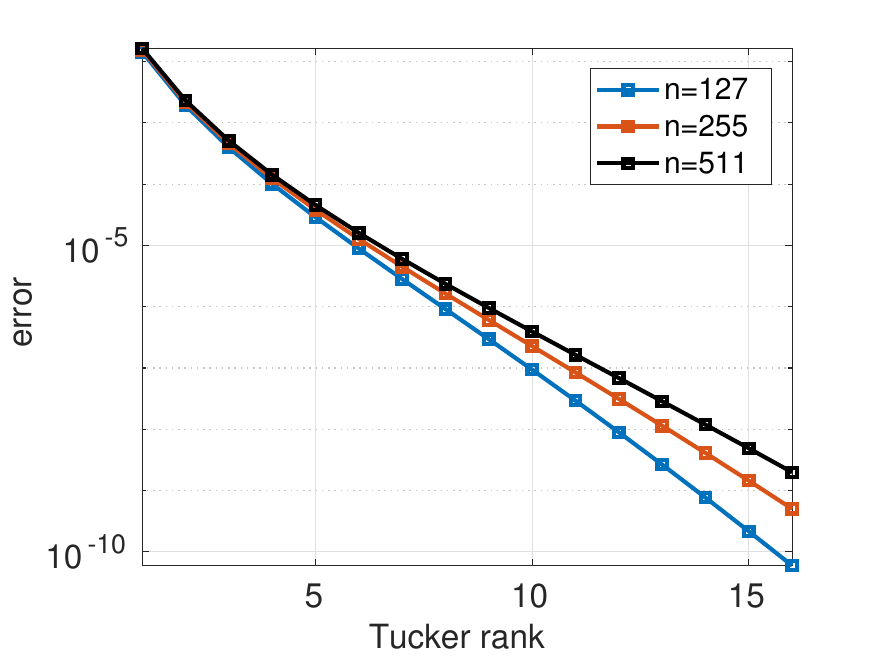}
\caption{\small Multigrid Tucker approximation of ${\bf G}_2$  for $\alpha=  1/2 $ (left)
and $\alpha=  1/10$ (right) vs. Tucker rank $r$ for   $n\times n \times n $ 
3D Cartesian grids with univariate grid size $n=127, \, 255, \,511$. }
\label{fig:Tucker_3D_MG}
\end{figure}

 The multigrid Tucker tensor decomposition for function related tensors \cite{KhKh3:08} 
 reduces the complexity 
 of the rank-${\bf r}$ Tucker tensor decomposition for full size tensors from $O(n^{d+1})$ to $O(n^{d})$. 
Here, the multigrid Tucker tensor approximation is used as a precomputing step 
for decomposing the 3D cores ${\bf G}_p=[g_p(i,j,k)]\in \mathbb{R}^{n\times n \times n}$, 
$p=1,\, 2,\, 3,$ in  (\ref{eq:coreTens}) and (\ref{eq:3a}) -- (\ref{eq:3c}),
corresponding to discretization on a sequence of $n_m\times n_m \times n_m $, $m=1,...,M$,
3D Cartesian grids. On the example of the 3D tensor ${\bf T}_0={\bf G}_{2,M}=[g_2(i,j,k)]$, 
defined by (\ref{eq:3b}),  Figure \ref{fig:Tucker_3D_MG} demonstrates  the exponentially fast decay
of the approximation error in the Tucker rank ${\bf r}$ (in Frobenius norm),
\[
 E_F=\frac{\| {\bf T}_0 -{\bf T}_{{\bf r},n}\|}{\|{\bf T}_0 \|},
 \]
for fractional powers $\alpha=  1/2$   (left) and for $\alpha=  1/10$ (right) calculated for 
$n=127,255,511$. 
Figure \ref{fig:Tucker_3D_MG} illustrates that for the 3D tensor ${\bf G}_2$
the separable representation with accuracy of the order of $\approx 10^{-7}$,
can be constructed using  rank-$10$ Tucker approximation, 
nearly independently on the size of discretization grid.  
The Tucker core of small size is transformed to a canonical tensor by 
the Tucker-to-canonical decomposition \cite{Khor2Book2:18},   
yielding the mixed tensor format, which is used as the starting rank-structured tensor
representation of the governing operator in the solution process. 
For reducing the Kronecker rank of the 
system matrix and of the current target vector in the course of the PCG iteration  
  the  RHOSVD-based  Tucker decomposition to the canonical tensors is applied \cite{Khor2Book2:18}.

\section{Numerics on rank-structured tensor numerical schemes}
\label{sec:numerics_tensor}

In this section we analyze the rank decomposition of all matrix entities 
involved in the solution operator (\ref{eqn:Lagrange_cont}). For the ease of exposition,
in what follows, we set the model constants as $\beta=\gamma=1$ and assume that
$n_1=n_2=n_3$.
Recall that $A=F^\ast G F$ with the notation $A=\mathcal{A}_h$, 
where $\mathcal{A}_h$ is the FDM approximation to the elliptic operator $\mathcal{A}$
and $G$ is the diagonal core matrix represented in terms of eigenvalues of the discrete 
Laplacian $A=\mathcal{A}_h$.
All numerical simulations are performed in MATLAB on a laptop.

\begin{figure}[tbh]
\centering
\includegraphics[width=5.0cm]{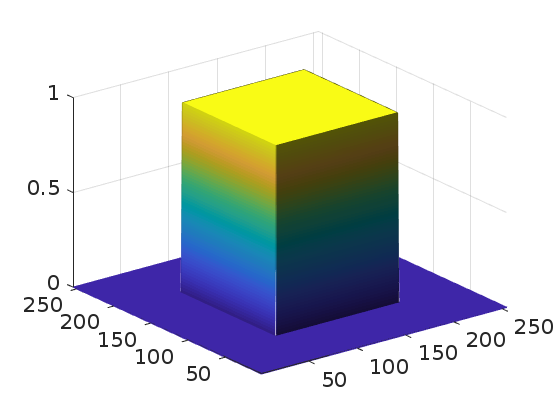}
\includegraphics[width=5.0cm]{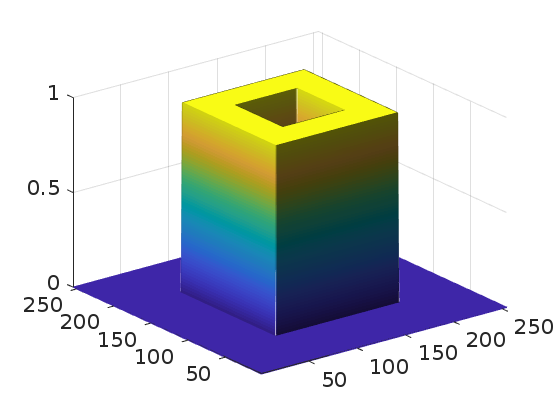}
\includegraphics[width=5.0cm]{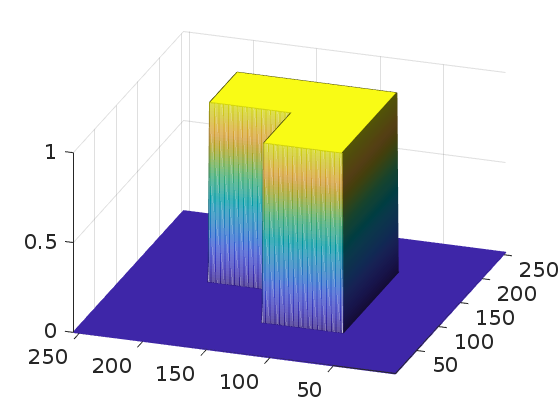}
\caption{\small Shapes of the right-hand side $y_\Omega$ used in the 2D 
equation (\ref{eqn:OptCont}) computed with n=255.}
\label{fig:design_RHS}
\end{figure}

\begin{figure}[tbh]
\centering
\includegraphics[width=5.0cm]{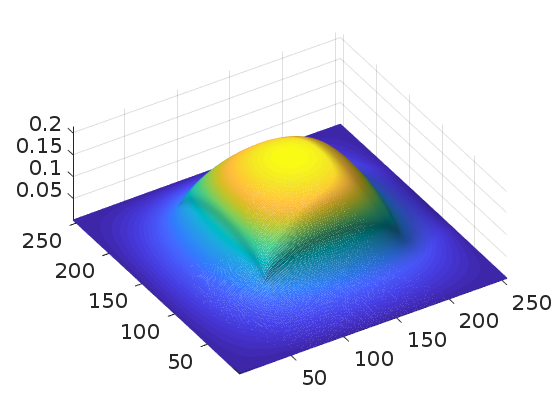}
\includegraphics[width=5.0cm]{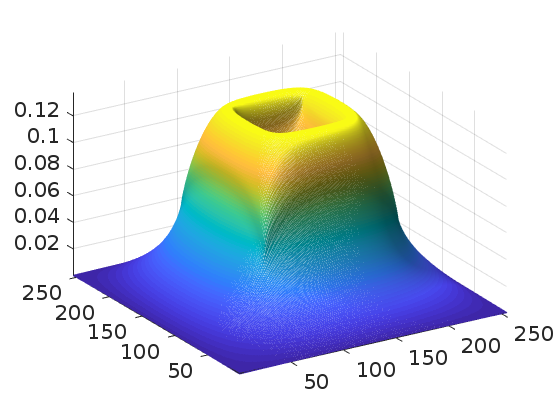}
\includegraphics[width=5.0cm]{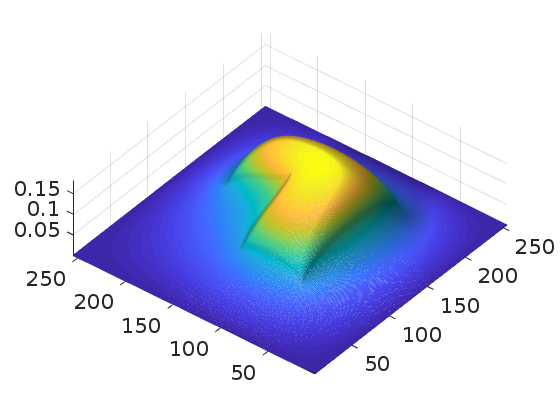}
\caption{\small Solutions ${\bf u}^\ast$ for above right-hand sides with $\alpha=1/2$ for $n=255$. }
\label{fig:FFT_2D_Alpha_12}
\end{figure}

\begin{figure}[tbh]
\centering
\includegraphics[width=5.0cm]{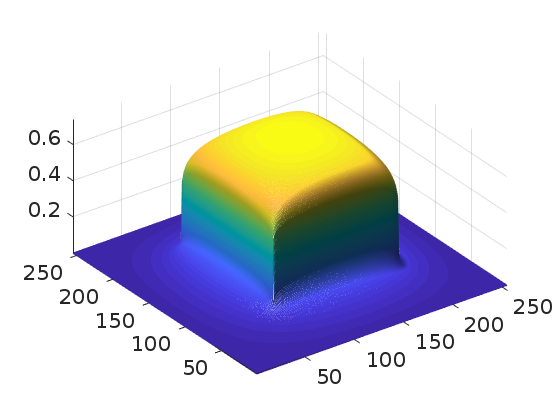}
\includegraphics[width=5.0cm]{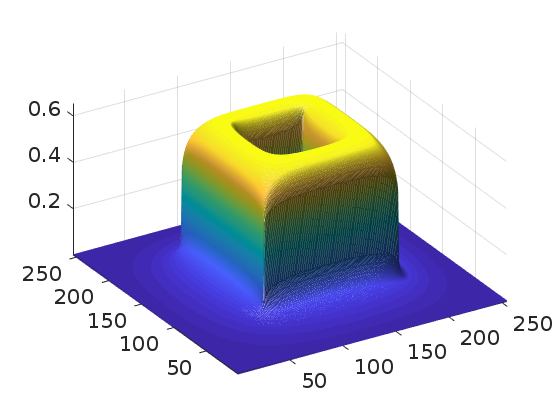}
\includegraphics[width=5.0cm]{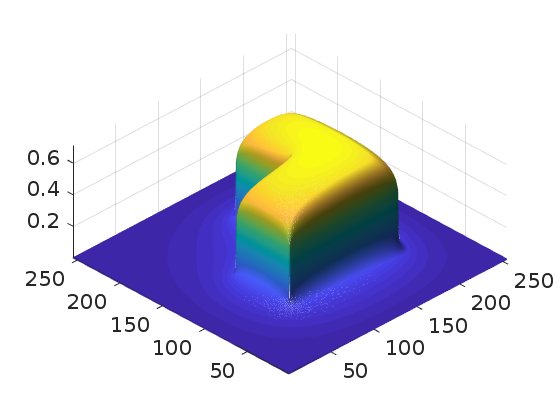}
\caption{\small Solutions ${\bf u}^\ast$ for above right-hand sides with $\alpha=1/10$ for $n=255$. }
\label{fig:FFT_2D_Alpha_110}
\end{figure}

 First, we illustrate the smoothing properties of the elliptic operator
 $\mathcal{A}_h^{-\alpha}$ in 2D
 (or by the other words, the localization properties of the fractional operator $\mathcal{A}_h^{\alpha}$)
 in the equation for control  depending on the fractional power $\alpha>0$.
 Figures \ref{fig:design_RHS}, \ref{fig:design_RHS}, \ref{fig:FFT_2D_Alpha_12}
 and \ref{fig:FFT_2D_Alpha_110} represent the shape of the design function $y_\Omega$
 and the corresponding optimal control ${\bf u}^\ast$ in the equation (\ref{eqn:OptCont})
 computed for different values of the parameter $\alpha$ and for fixed grid size $n=255$. 
   \begin{figure}[tbh]
\centering
\includegraphics[width=5.0cm]{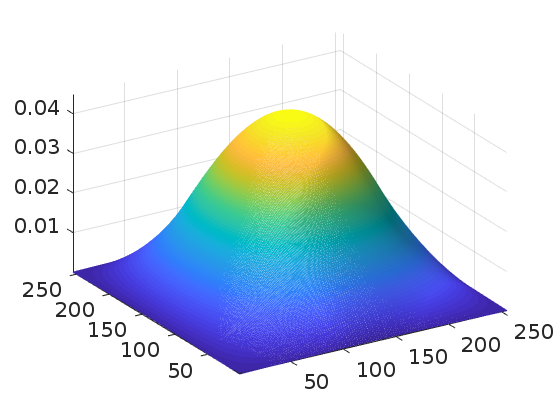}
\includegraphics[width=5.0cm]{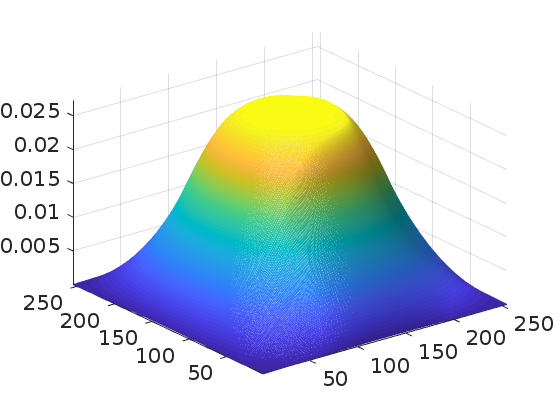}
\includegraphics[width=5.0cm]{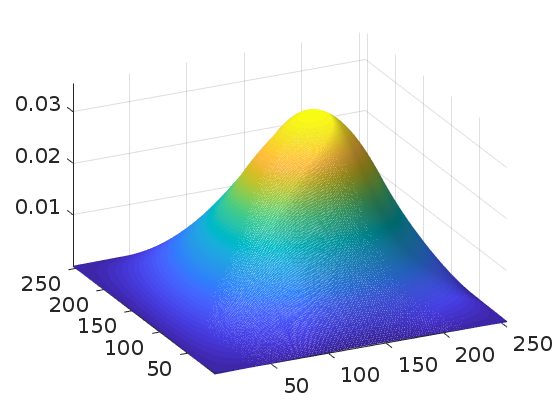}
\caption{\small Solutions ${\bf u}^\ast$ for above right-hand sides $y_\Omega$ 
with $\alpha=1$ for $n=255$. }
\label{fig:FFT_2D_Alpha_1}
\end{figure}

One observes the nonlocal features of the 
 elliptic resolvent $\mathcal{A}_h^{-1}$ and highly localized action of the operator
 $\mathcal{A}_h^{-\frac{1}{10}}$.

\subsection{Numerical tests for 2D case}
\label{ssec:numerics_tensor_2D}

\begin{figure}[tbh]
\centering
\includegraphics[width=6.4cm]{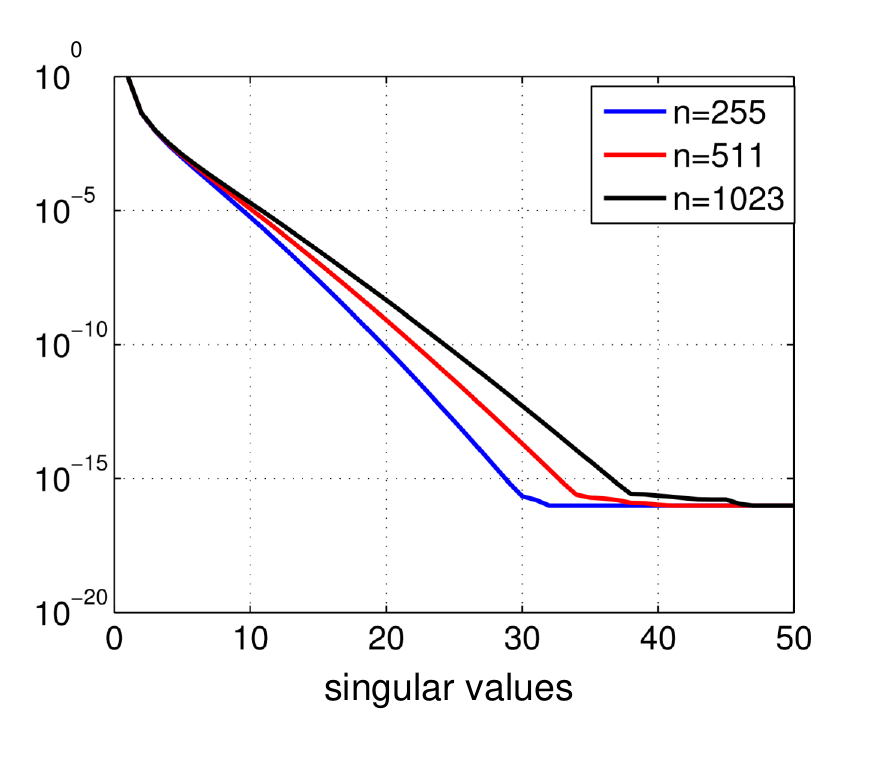}
\includegraphics[width=6.4cm]{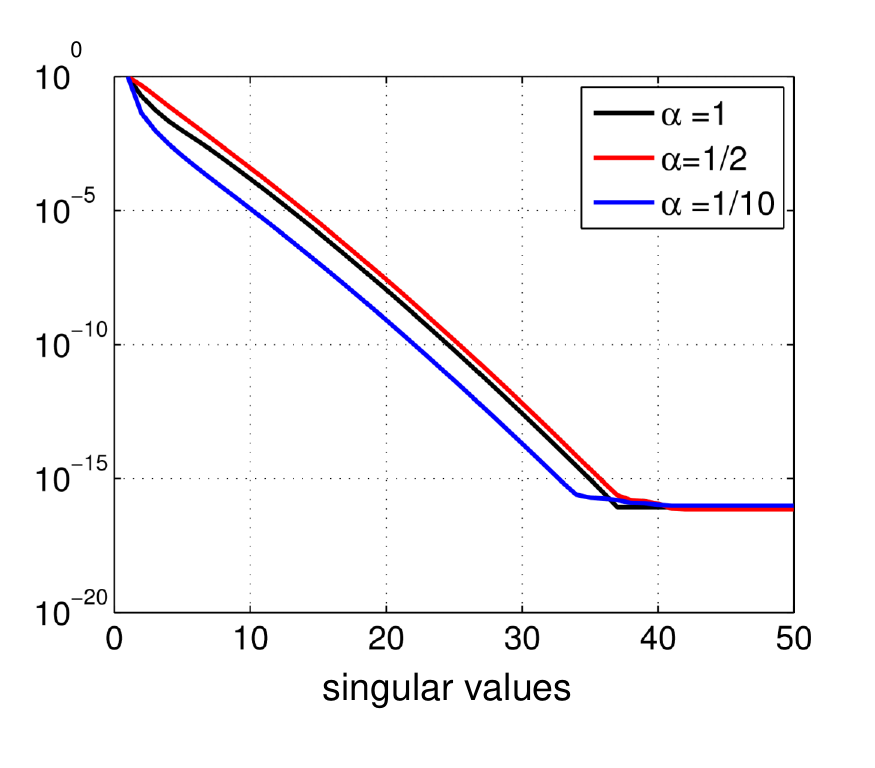}
\caption{Decay of singular values for $G_1$ with $\alpha=1$ vs. $n$ (left); 
singular values for $G_1$ vs. $\alpha>0$ with fixed $n=511$ (right). }
\label{fig:SVD_2D_Grid_Alpha}
\end{figure}

 Figure \ref{fig:SVD_2D_Grid_Alpha}, left, represents the singular values of the matrix 
 $G_1$, with entries given by (\ref{eq:2a}) for different univariate grid size 
 $n=255,511$, and $1023$ and fixed $\alpha=1$ (Laplacian inverse).
 Figure \ref{fig:SVD_2D_Grid_Alpha}, right, shows the decay of respective singular values
 for $G_1$ with fixed univariate grid size $n=511$ and for different $\alpha=1, 1/2,1/10$.

\begin{figure}[tbh]
\centering
\includegraphics[width=6.4cm]{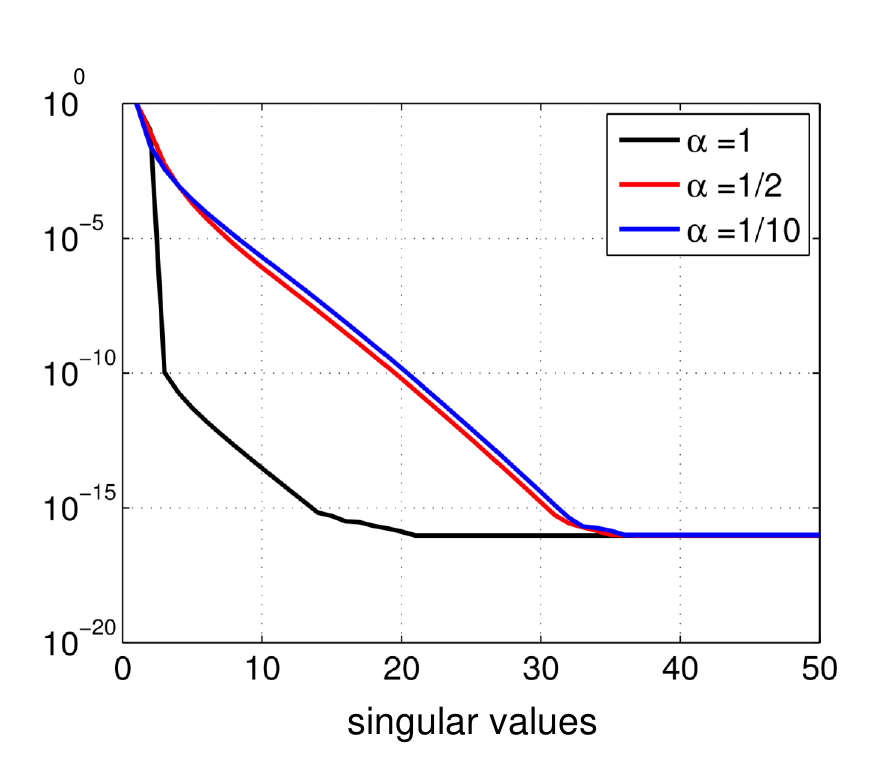}
\includegraphics[width=6.4cm]{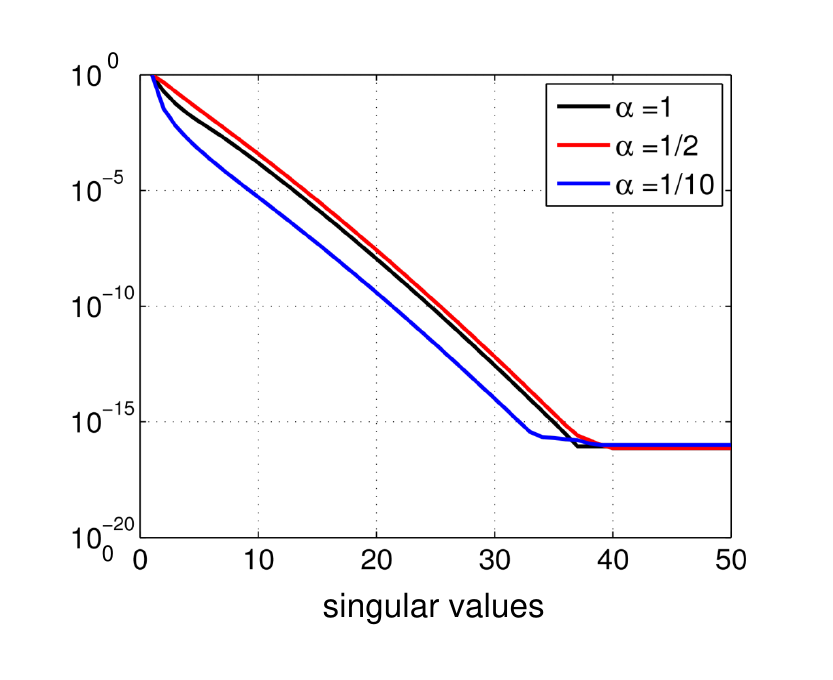}
\caption{Decay of singular values of $G_2$ (left) and $G_3$ (right) vs. 
$\alpha=1, 1/2,1/10$ for $n=511$. }
\label{fig:SVD_2D_Alpha}
\end{figure}

Figure \ref{fig:SVD_2D_Alpha} demonstrates the behavior of singular values for matrices
 $G_2$ and $G_3$, with the entries corresponding to (\ref{eq:2b}) and (\ref{eq:2c}), respectively, 
 vs. $\alpha=1, 1/2,1/10$ with fixed univariate grids size $n=511$.
 In all cases we observe exponentially fast decay of the singular values which means 
 there exists the accurate low Kronecker rank approximation of the matrix functions $A_1, A_2$ and $A_3$
 (see equations (\ref{eqn:L_m_alpha}), (\ref{eqn:L_alp_m_alp}) and (\ref{eqn:(L_alp_m_alp)_m1}))
 including fractional powers of the elliptic operator.

\begin{figure}[tbh]
\centering
\includegraphics[width=6.4cm]{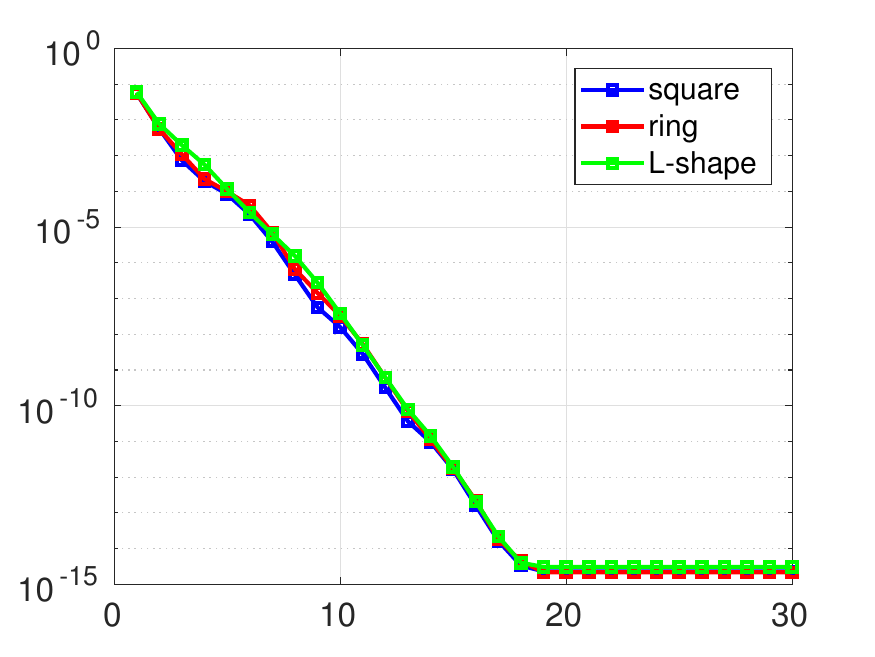}
\caption{\small Decay of the error for the optimal control vs. truncation rank parameter.}
\label{fig:OpContsol_rank_decay}
\end{figure}

Decay of the error for the optimal control obtained 
as the solution of equation (\ref{eqn:OptCont}) with rank-$R$ approximation of the solution operator
$A_3$ is shown in Figure \ref{fig:OpContsol_rank_decay}.

\begin{figure}[tbh]
\centering
\includegraphics[width=10.0cm]{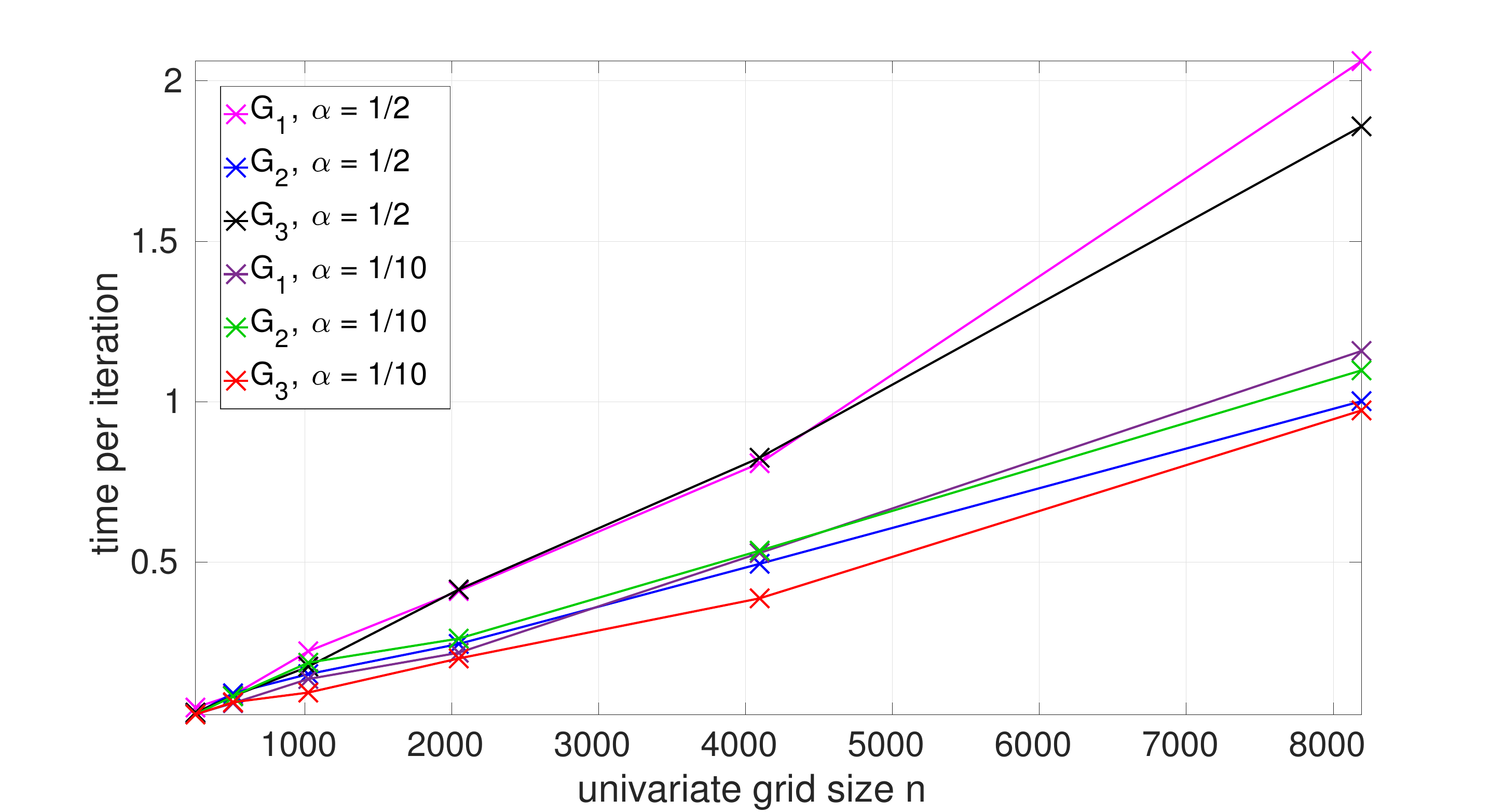}
\caption{\small CPU times (sec) vs. univariate grid size $n$ for a single iteration 
of Algorithm~\ref{alg:pcg} for a 2D problem, for different fractional operators 
and fixed preconditioner rank $r = 5$.}
\label{fig:times2D}
\end{figure}

As we have shown theoretically in Section~\ref{sec:Larrange_eqn}, a single PCG iteration has a
complexity, which is slightly higher than linear in the univariate grid size $n$. Figure \ref{fig:times2D}
shows that the CPU times show the expected behavior. Thus, with Figure~\ref{fig:times2D} and
Tables~\ref{tab:precond2D_alp12}~and~\ref{tab:precond2D_alp110}, the overall cost of the algorithm
is almost linear in the univariate grid size $n$ for the problem discretized 
on $n \times n$ 2D Cartesian grid.

We also test the properties of the low-rank discrete operator as a preconditioner. This means,
we solve the equations in $\mathbb{R}^d$, $d=2,\, 3$,
\begin{align}
	A^{\alpha} \mathbf{x} &= \mathbf{b},\label{eq:precond_ex1}\\
	\big( I + A^{2\alpha} \big) \mathbf{x} &= \mathbf{b},\label{eq:precond_ex2}\\
	\big( A^{\alpha} + A^{-\alpha} \big) \mathbf{x} &= \mathbf{b},\label{eq:precond_ex3}
\end{align}
with a preconditioned conjugate gradient scheme, 
using a low-rank direct solver as a preconditioner  discussed above. We simplify the notation
by $A={\cal A}_h$.

In numerical tests we solve the equations \eqref{eq:precond_ex1} - \eqref{eq:precond_ex3}
on a grid of size $n$, using a rank-$r$ preconditioner.
Tables \ref{tab:precond2D_alp12} and \ref{tab:precond2D_alp110} show the number of 
CG iteration counts for convergence to a relative residual of $10^{-6}$ 
of \eqref{eq:precond_ex1}-\eqref{eq:precond_ex3} with $\alpha = 1/2$ and
$\alpha = 1/10$, respectively. The dash `---' indicates failure to converge to 
converge in 100 iterations.
\begin{table}[tbh]
\begin{center}\footnotesize
\begin{tabular}
  [c]{|r|r|r|r|r||r|r|r|r||r|r|r|r|r||}%
\hline
 & \multicolumn{4}{|c|}{$g_1$} & \multicolumn{4}{|c|}{$g_4$ } & \multicolumn{4}{|c|}{$g_3$} \\
\hline 
\diagbox{$r$}{$n$}  & 256 & 512 & 1024 & 2048 & 256  & 512   & 1024  & 2048  &  256  & 512 & 1024 & 2048  \\  
\hline  
1                 & 20 & 24 & 24  & 29  & ---   & ---   & 83   & 80 & 20  & 24  & 24 & 19 \\
\hline
2                 & --- & ---  & 3 & 3 & 73  & ---   & 38  & 36 & ---  & ---  & 3 & 3\\
\hline
3                  & 7 & 9  & 10  & 14 & 99  & ---  & 17  & 16  & 7  & 9  & 10 & 14 \\
\hline
4                  & 5  & 6  & 6  & 9 & 31  & ---  & 3  & 3  & 5   & 6 & 6 & 9 \\
\hline
5                  & 4  &  4 &  4 & 5  & 11  & --  & 5  & 5 & 4   & 4  & 4 & 5 \\
\hline
6                  & 3  &  3  & 3 & 4  & 6   & 13  & 2  & 2  & 3   & 3 & 3 & 4\\
\hline
7                  & 3  &  3  & 3 & 3  & 4   & 7   & 6  & 4  & 3   & 3 & 3 & 3\\
\hline
8                  & 2  &  2  & 2 & 2  & 3   & 5   & 4  & 2  & 2   & 2 & 2 & 2\\
\hline
9                  & 2  &  2  & 2 & 2  & 3   & 4   & 3  & 4  & 2   & 2 & 2 & 2 \\
\hline
10                 & 2  &  2  & 2 & 2  & 3   & 3   & 2  & 3  & 2   & 2  & 2 & 2\\
\hline
 \end{tabular}
 \caption{\small PCG iteration counts for convergence to a relative residual 
 of $10^{-6}$ for the equations \eqref{eq:precond_ex1} - \eqref{eq:precond_ex3}
 for a 2D fractional Laplacian with $\alpha = 1/2$ 
 vs. the univariate grid size $n$ and separation rank $r$.}
   \label{tab:precond2D_alp12}
\end{center}
\end{table}
   
As can be seen in Tables~\ref{tab:precond2D_alp12}~and~\ref{tab:precond2D_alp110}, we achieve almost
grid-independent preconditioning; the iteration counts only grow logarithmically with the number of grid 
points, as can be expected from the theoretical reasoning. As can be seen in Table~
\ref{tab:precond2D_alp12}, the ranks of the preconditioner should be chosen sufficiently large to ensure
reliability. In the cases tested here, $r=6$ is sufficient to achieve reliable preconditioning even in
the most difficult case of equation \eqref{eq:precond_ex2} with $\alpha=1/2$.
  \begin{table}[tbh]
\begin{center}\footnotesize
\begin{tabular}
  [c]{|r|r|r|r|r||r|r|r|r||r|r|r|r|r||}%
\hline
 & \multicolumn{4}{|c|}{$g_1$} & \multicolumn{4}{|c|}{$g_4$ } & \multicolumn{4}{|c|}{$g_3$} \\
\hline 
\diagbox{$r$}{$n$} & 256 & 512 & 1024 & 2048 & 256 & 512 & 1024 & 2048 & 256 & 512 & 1024 & 2048 \\  
\hline  
1                  & 9 & 9 & 10 & 11 & 11  & 13  & 14  & 16  & 7  & 7  & 8 & 9\\
\hline
2                  & 6 & 4 & 7 & 8 &   7   & 8   & 8  & 9  & 5  & 5  & 6 & 6\\
\hline
3                  & 4 & 5 & 5 & 6 &   5   & 5   & 6  & 7  & 4  & 4  & 5 & 5\\
\hline
4                  & 4 & 4 & 4 & 5 &   4   & 4   & 4  & 5  & 3  & 4 & 4 & 4\\
\hline
5                  & 3 & 3 & 4 & 4  &  3   & 4   & 4  & 4 & 3   & 3  & 3 & 4\\
\hline
6                  & 3 & 3 & 3 & 4  &  3   & 3   & 3  & 4  & 2   & 3 & 3 & 3\\
\hline
7                  & 2 & 3 & 3 & 3 &   2   & 3   & 3  & 3  & 2   & 2 & 3 & 3\\
\hline
8                  & 2 & 2 & 2 & 3 &   2   & 2   & 2  & 3  & 2   & 2 & 2 & 3 \\
\hline
9                  & 2 & 2 & 2 & 2  &  2   & 2   & 2  & 2  & 2   & 2 & 2 & 3\\
\hline
10                 & 2 & 2 & 2 & 2  &  2   & 2   & 2  & 2  & 2   & 2  & 2 & 2\\
\hline
 \end{tabular}
 \caption{\small PCG iteration counts for convergence to a relative residual 
 of $10^{-6}$ for the equations \eqref{eq:precond_ex1} - \eqref{eq:precond_ex3}
 for a 2D fractional Laplacian with $\alpha = 1/10$ vs. 
 the univariate grid size $n$ and separation rank $r$.}
    \label{tab:precond2D_alp110}
\end{center}
\end{table}

 \subsection{Numerical tests for 3D case}
\label{ssec:numerics_tensor_3D}

  In the following examples we solve the problems governed by the 3D operators in  
  (\ref{eqn:L_m_alpha}) -- (\ref{eqn:(L_alp_m_alp)_m1}),   
  with a 3D fractional Laplacian with $\alpha=1, 1/2$ and $\alpha =1/10$.
  The rank-structured approximation to the above fractional operators is performed by using
  the multigrid Tucker decomposition of the 3D tensors ${\bf G}_k$, $k=1,2,3,4$,
  described by (\ref{eq:3a}) -- (\ref{eq:3c}), and the consequent Tucker-to-canonical
  decomposition of the Tucker core tensor thus obtaining a canonical tensor with
  a smaller rank. The rank truncation procedure in the PCG Algorithm~\ref{alg:pcg} is performed by
  using the RHOSVD tensor approximation and its consequent transform to the canonical format, see 
  Section \ref{sec:low-rank_tensor}.
  
\begin{figure}[tbh]
\centering
\includegraphics[width=7.0cm]{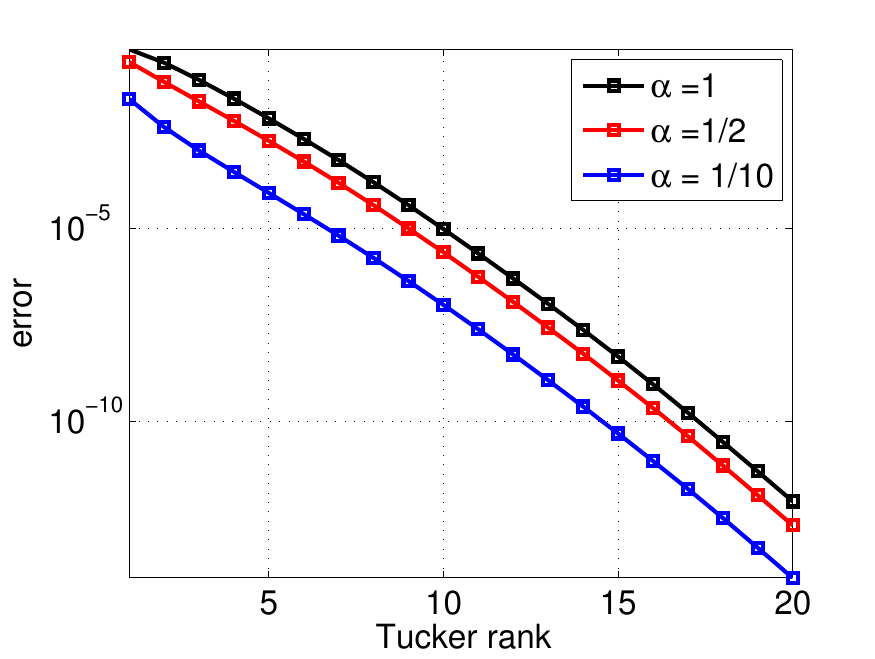}
\caption{\small Tucker tensor approximation of ${\bf G}_1$ vs. rank parameter for
$\alpha=1, 1/2,1/10$. }
\label{fig:Tucker_3D_Alpha_Inv}
\end{figure}

\begin{figure}[tbh]
\centering
\begin{tabular}{cc}
	\includegraphics[width=6.5cm]{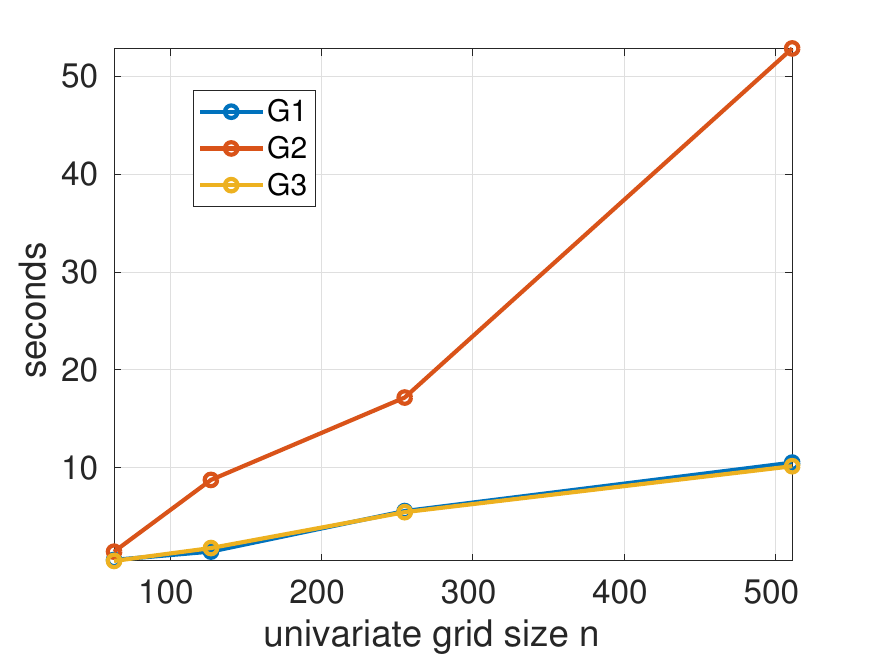}&
	\includegraphics[width=6.5cm]{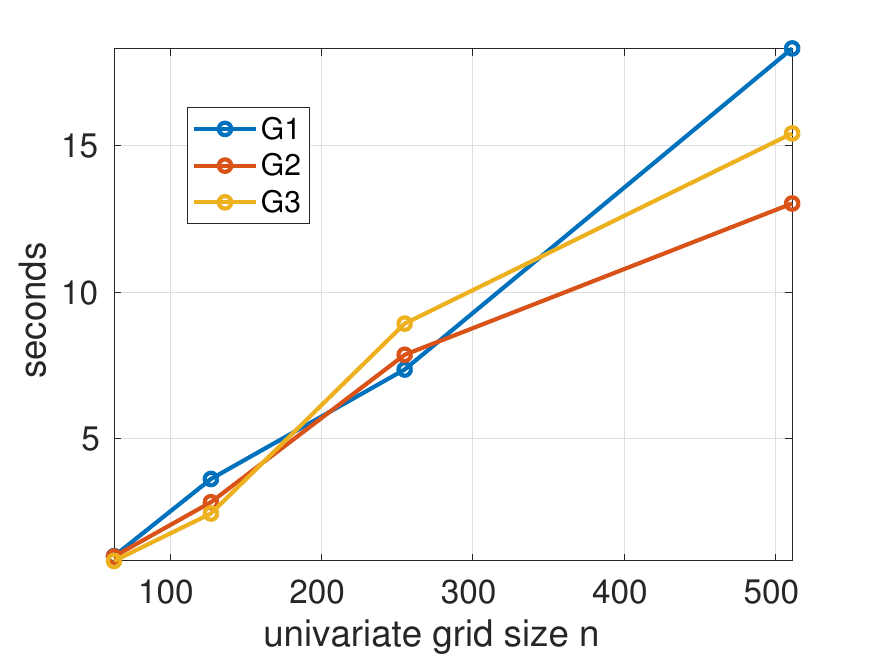}\\
	$\alpha = 1/2$, $r=4$&
	$\alpha = 1/10$, $r=7$
\end{tabular}
\caption{\small CPU times (in seconds) vs. univariate grid size $n$ for a single iteration of Algorithm~\ref{alg:pcg} 
for a 3D problem, for different fractional operators and fixed preconditioner rank $r$.}
\label{fig:times3D}
\end{figure}

\begin{figure}[tbh]
\centering
\includegraphics[width=7.0cm]{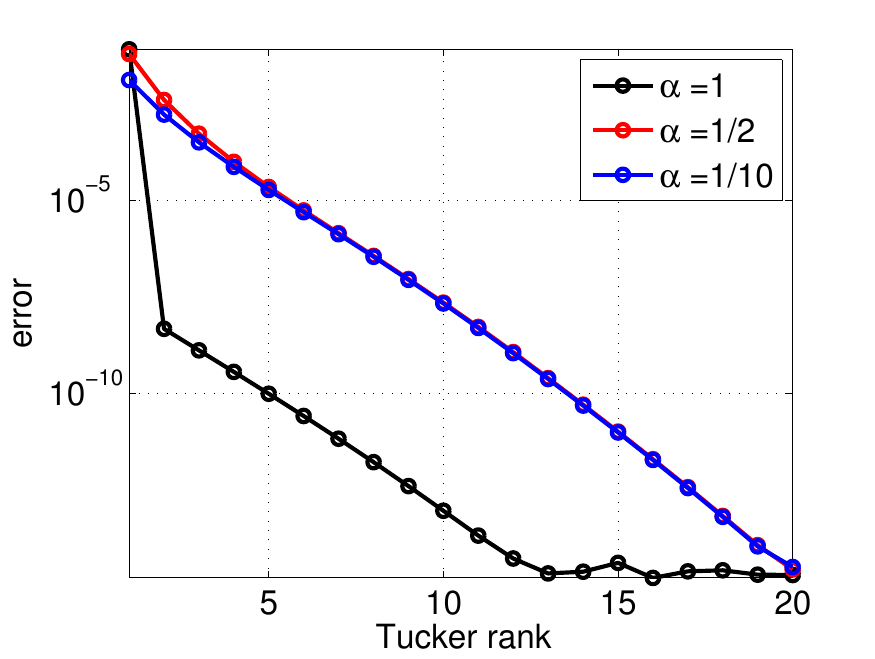}
\includegraphics[width=7.0cm]{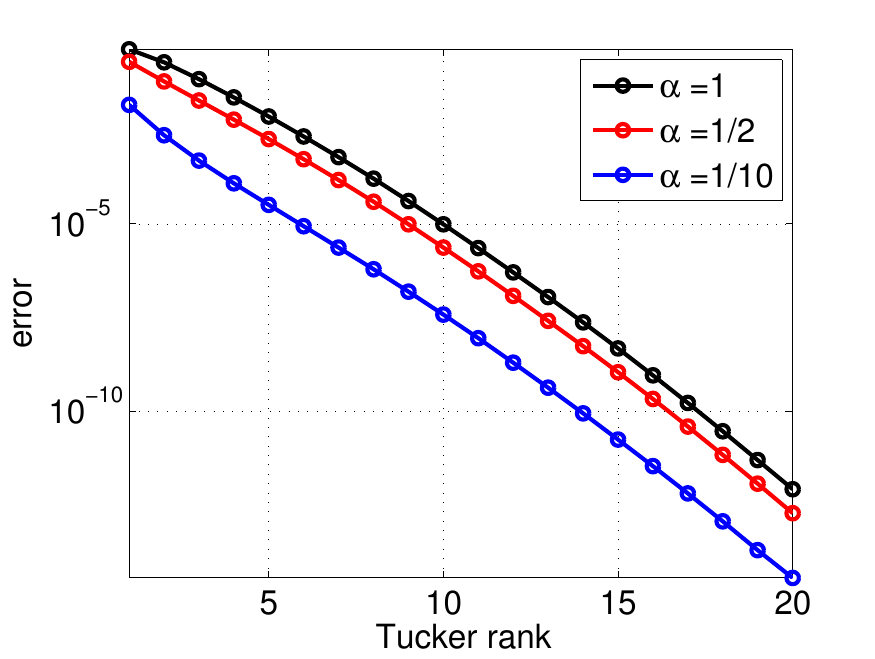}
\caption{\small Tucker tensor approximation of ${\bf G}_2$ and ${\bf G}_3$ vs. rank parameter 
for $\alpha=1, 1/2,1/10$. }
\label{fig:Tucker_3D_Alpha_SolOp}
\end{figure}

  Figures \ref{fig:Tucker_3D_Alpha_Inv} -- \ref{fig:Tucker_3D_Alpha_SolOp} demonstrate the
  exponential convergence of the approximation error with respect to the Tucker rank
  for operators given by \eqref{eq:precond_ex1} -- \eqref{eq:precond_ex3}.

We solve the equations \eqref{eq:precond_ex1} - \eqref{eq:precond_ex3}
using $ n\times n \times n$ 3D Cartesian 
grids with the univariate grid size $n$, using a rank-$r$ preconditioner.
Tables \ref{tab:precond3D_alp12} and \ref{tab:precond3D_alp110} show the number of 
CG iteration counts for convergence to a relative residual of $10^{-6}$ 
of \eqref{eq:precond_ex1} - \eqref{eq:precond_ex3} with $\alpha = 1/2$ and
$\alpha = 1/10$, respectively.

Similarly to the previous subsection, we see that the low-rank approximation gives us an approximately
grid-independent preconditioner. In the cases tested here, $r=6$ is sufficient to achieve reliable 
preconditioning even in the most difficult case of equation \eqref{eq:precond_ex2} with $\alpha=1/2$.

\begin{table}[tbh]
\begin{center}\footnotesize
\begin{tabular}
  [c]{|r|r|r|r|r||r|r|r|r||r|r|r|r|r||}%
\hline
 & \multicolumn{4}{|c|}{$g_1$} & \multicolumn{4}{|c|}{$g_4$ } & \multicolumn{4}{|c|}{$g_3$} \\
\hline 
\diagbox{$r$}{$n$}  & 64 & 128& 256 & 512 & 64 & 128 & 256 & 512 & 64 & 128 & 256 & 512  \\  
\hline
4                   & 1  & 2  & 1   & 1   & 1  & 6   & 1   & 2   & 1  & 2   & 1   & 1 \\
\hline
5                   & 1  & 1  & 1   & 2   & 1  & 1   & 8   & 4   & 1  & 1   & 1   & 2 \\
\hline
6                   & 1  &  1 & 1   & 1   & 2  & 2   & 1   & 1   & 1  & 1   & 1   & 1\\
\hline
7                   & 1  &  3 & 1   & 2   & 1  & 1   & 5   & 4   & 1  & 2   & 1   & 2\\
\hline 
8                   & 1  &  1 & 1   & 1   & 1  & 1   & 1   & 1   & 1  & 1   & 1   & 1\\
\hline
9                   & 1  &  1 & 1   & 2   & 1  & 6   & 5   & 4   & 1  & 1   & 1   & 2 \\
\hline
10                  & 1  &  1 & 1   & 1   & 1  & 6   & 1   & 1   & 1  & 1   & 1   & 1\\
\hline
 \end{tabular}
 \caption{\small PCG iteration counts for convergence to a relative residual 
 of $10^{-6}$ for the equations \eqref{eq:precond_ex1} - \eqref{eq:precond_ex3}
 for a 3D fractional Laplacian with $\alpha = 1/2$.
 Here $n$ is the univariate grid size, $r$ is the separation rank.}
   \label{tab:precond3D_alp12}
\end{center}
\end{table}

  \begin{table}[tbh]
\begin{center}\footnotesize
\begin{tabular}
  [c]{|r|r|r|r|r||r|r|r|r||r|r|r|r|r||}%
\hline
 & \multicolumn{4}{|c|}{$g_1$} & \multicolumn{4}{|c|}{$g_4$ } & \multicolumn{4}{|c|}{$g_3$} \\
\hline 
\diagbox{$r$}{$n$}  & 64 & 128& 256 & 512 & 64 & 128 & 256 & 512 & 64 & 128 & 256 & 512  \\  
\hline
4                   & 2  & 1  & 9   & 20  & 2  & 1   & 10  & 17  & 1  & 1   & 9   & 18 \\
\hline
5                   & 1  & 1  & 1   & 1   & 1  & 1   & 1   & 1   & 1  & 2   & 1   & 13 \\
\hline
6                   & 1  &  1 & 1   & 2   & 1  & 1   & 1   & 2   & 1  & 1   & 1   & 7 \\
\hline
7                   & 1  &  1 & 1   & 2   & 1  & 1   & 1   & 2   & 1  & 1   & 2   & 1 \\
\hline 
8                   & 1  &  1 & 1   & 1   & 1  & 1   & 1   & 1   & 1  & 1   & 1   & 2 \\
\hline
9                   & 1  &  1 & 1   & 1   & 1  & 1   & 1   & 1   & 1  & 1   & 1   & 1 \\
\hline
10                  & 1  &  1 & 1   & 2   & 1  & 1   & 1   & 2   & 1  & 1   & 1   & 1 \\
\hline
 \end{tabular}
 \caption{\small CG iteration counts for convergence to a relative residual 
 of $10^{-6}$ for the equations \eqref{eq:precond_ex1}-\eqref{eq:precond_ex3}
 for a 3D fractional Laplacian with $\alpha = 1/10$.
 Here $n$ is the univariate grid size, $r$ is the separation rank.}
    \label{tab:precond3D_alp110}
\end{center}
\end{table}
  
 Our numerical test indicates that all three matrices  $A_1, A_2$ and $A_3$, as well as 
 the corresponding three-tensors have 
 $\varepsilon$-rank approximation such that the rank parameter depends logarithmically on
 $\varepsilon$, i.e., $r=O(|\log \varepsilon|)$, that
 means that the low rank representation of the design function 
 $y_\Omega$ ensures the low rank representation of both optimal control and optimal 
 state variable. 
 
 We show as well that, using rank-structured tensor methods for the numerical solution 
 of this optimization problem using the operators of type $A_1, A_2$ and $A_3$  
 can be implemented at low cost that scales linearly in the univariate grid size, $O(n \log n)$, see 
 Figure~\ref{fig:times3D}.
 
   \begin{figure}[h]
\centering
\includegraphics[width=5.4cm]{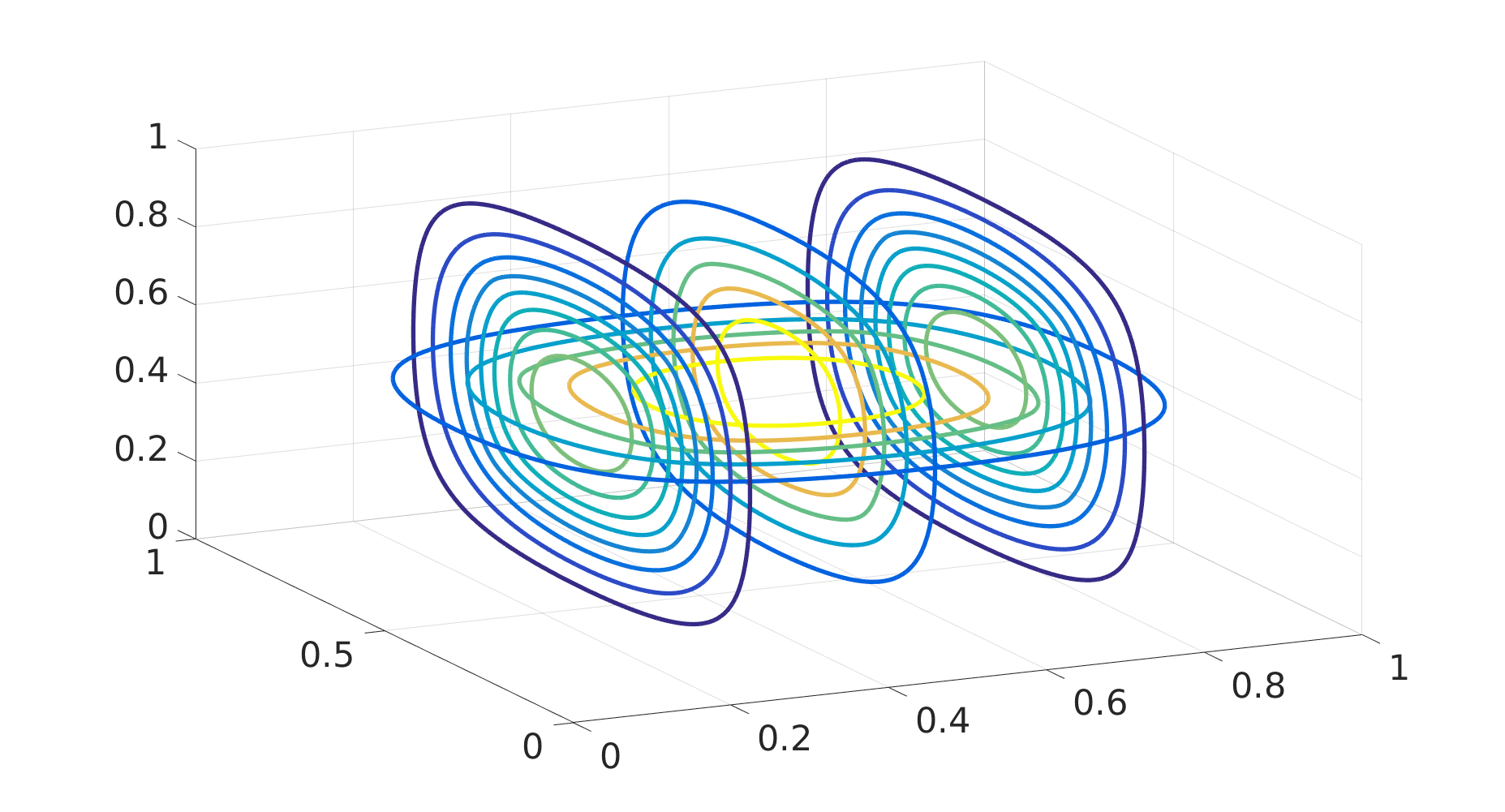}
\includegraphics[width=5.4cm]{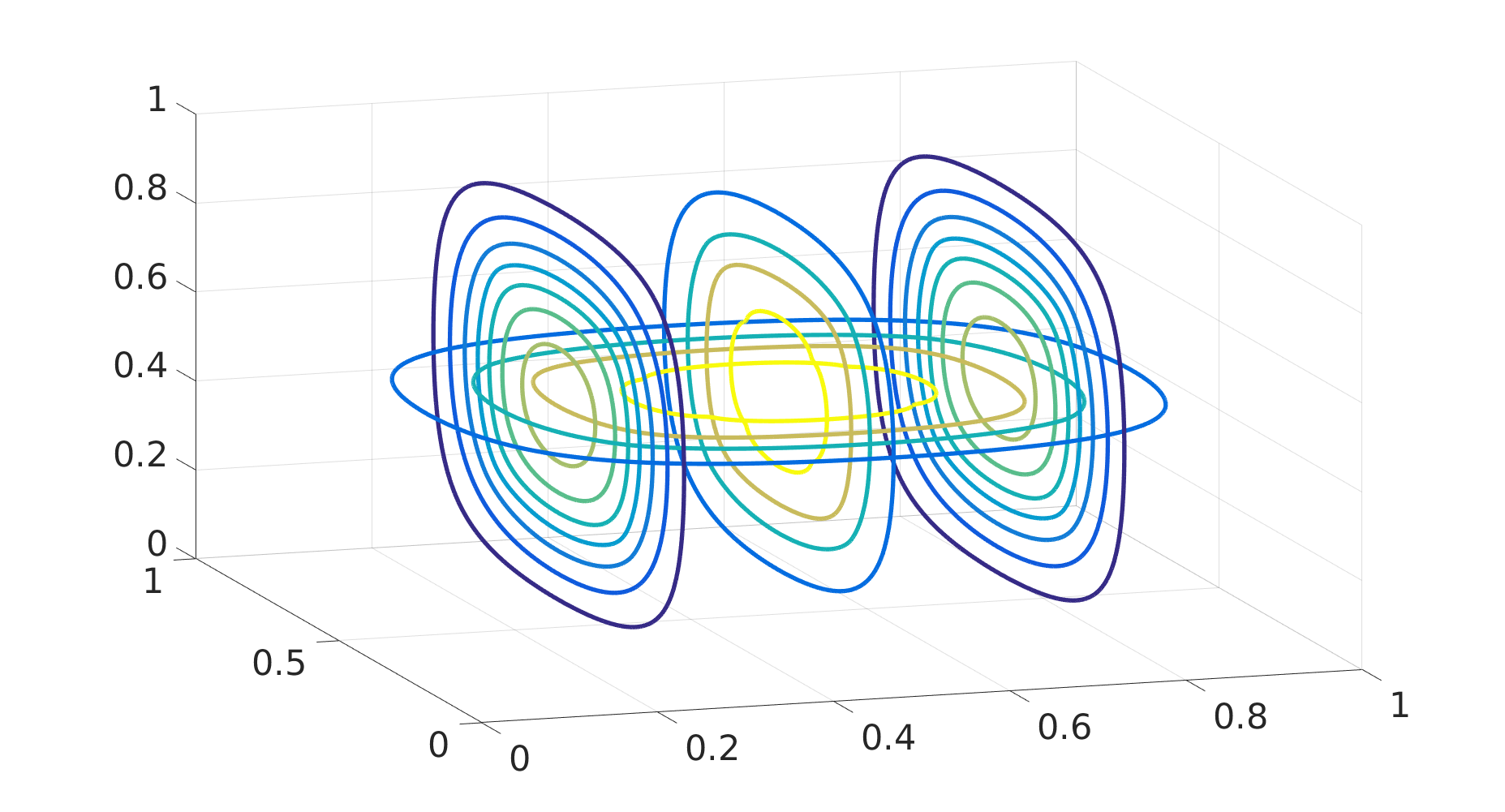}
\includegraphics[width=5.4cm]{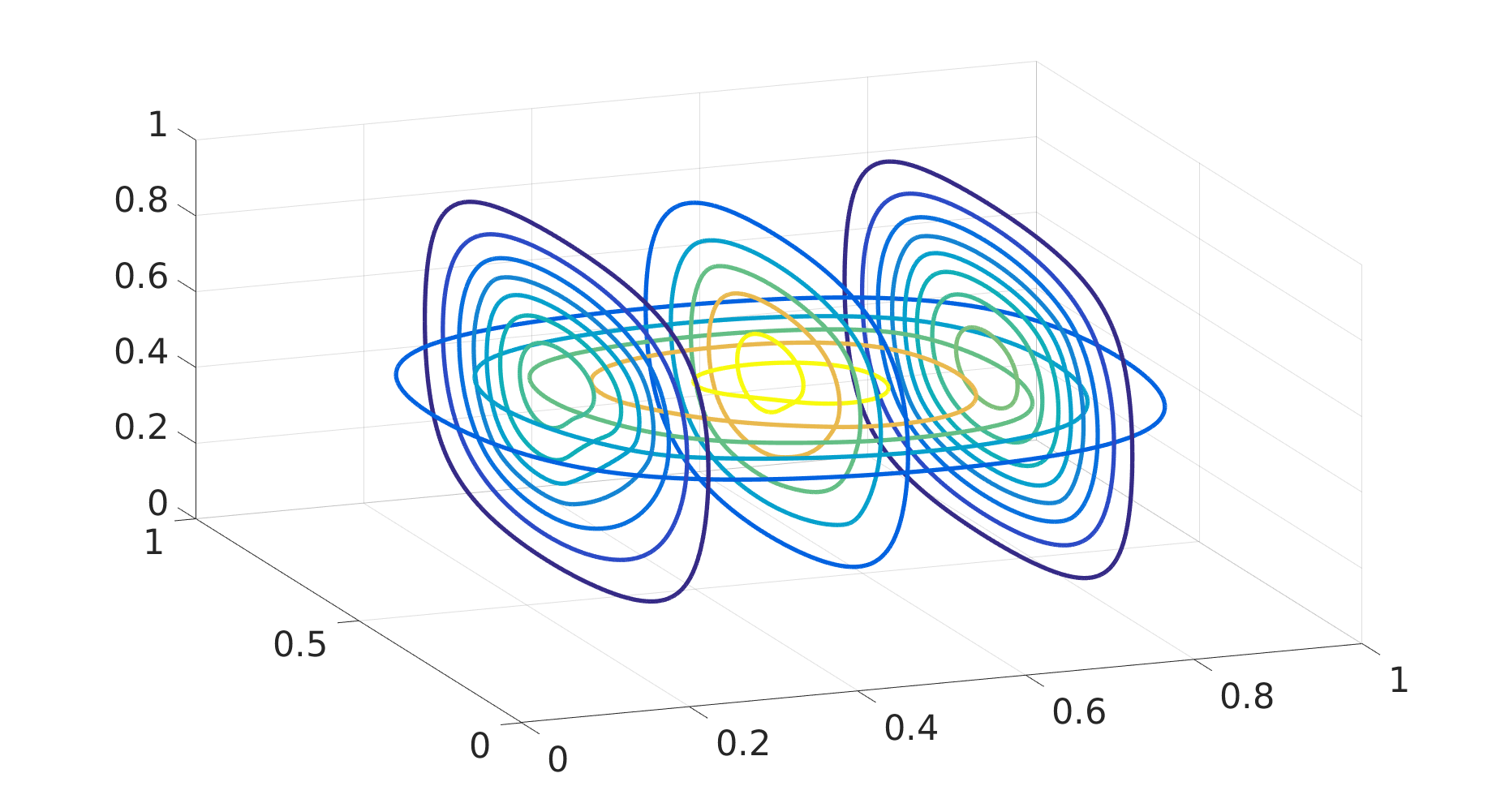}
\caption{\small Solutions of the equation with 3D right-hand sides (analogous to
Figure \ref{fig:design_RHS}) with $\alpha=1$ for $n=255$.}
\label{fig:FFT_3D_Alpha_1}
\end{figure}

   \begin{figure}[h]
\centering
\includegraphics[width=5.4cm]{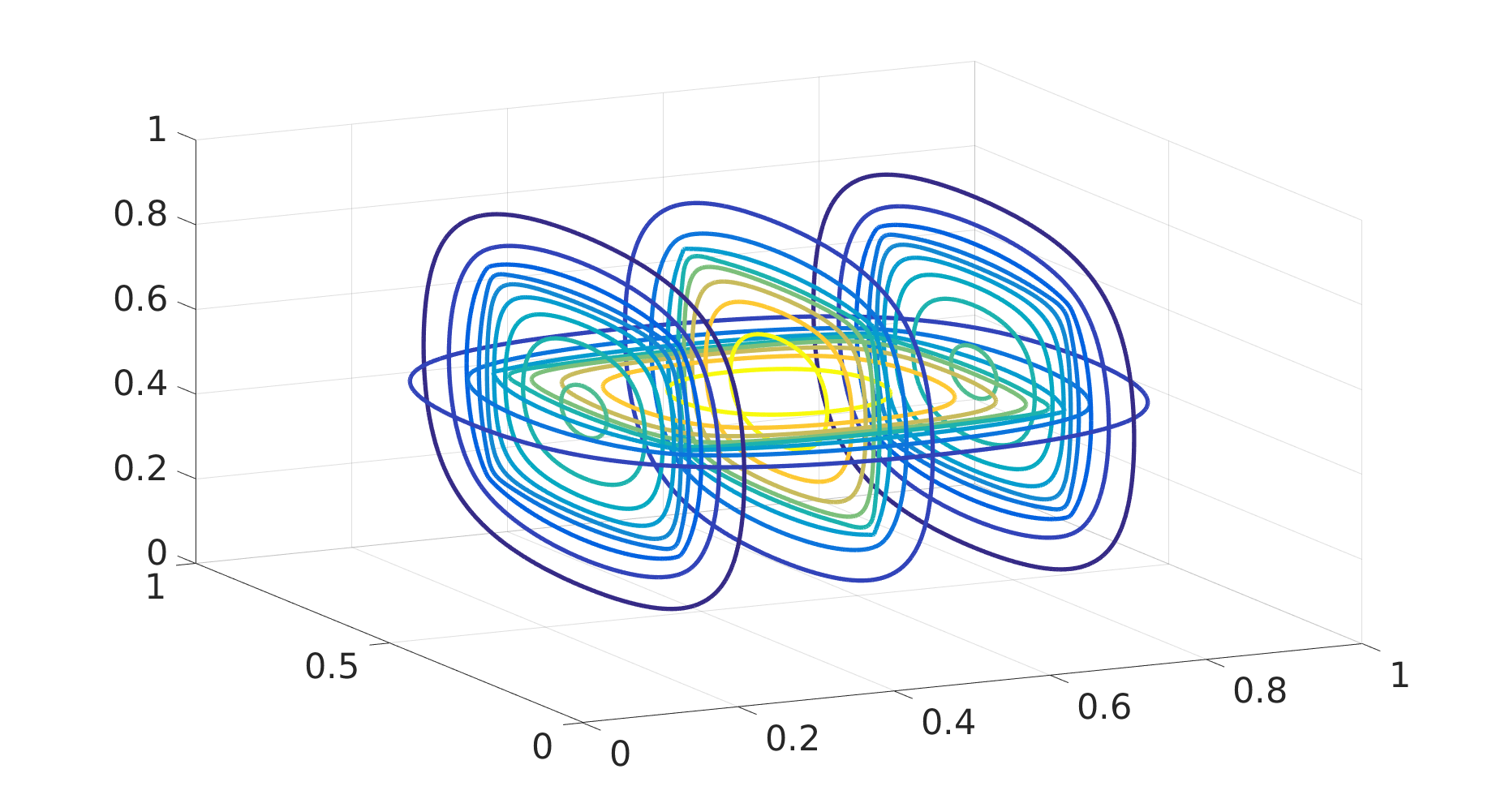}
\includegraphics[width=5.4cm]{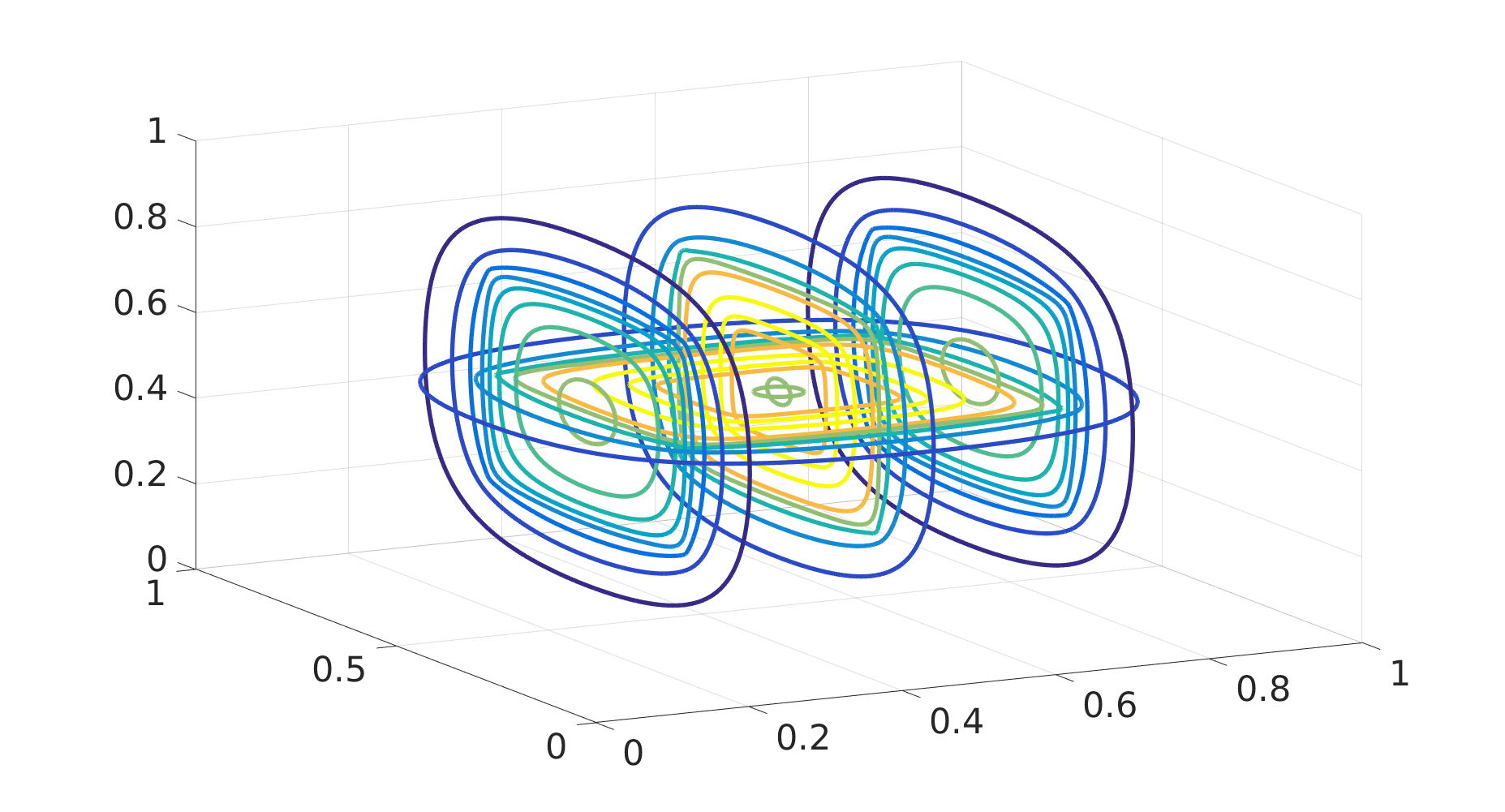}
\includegraphics[width=5.4cm]{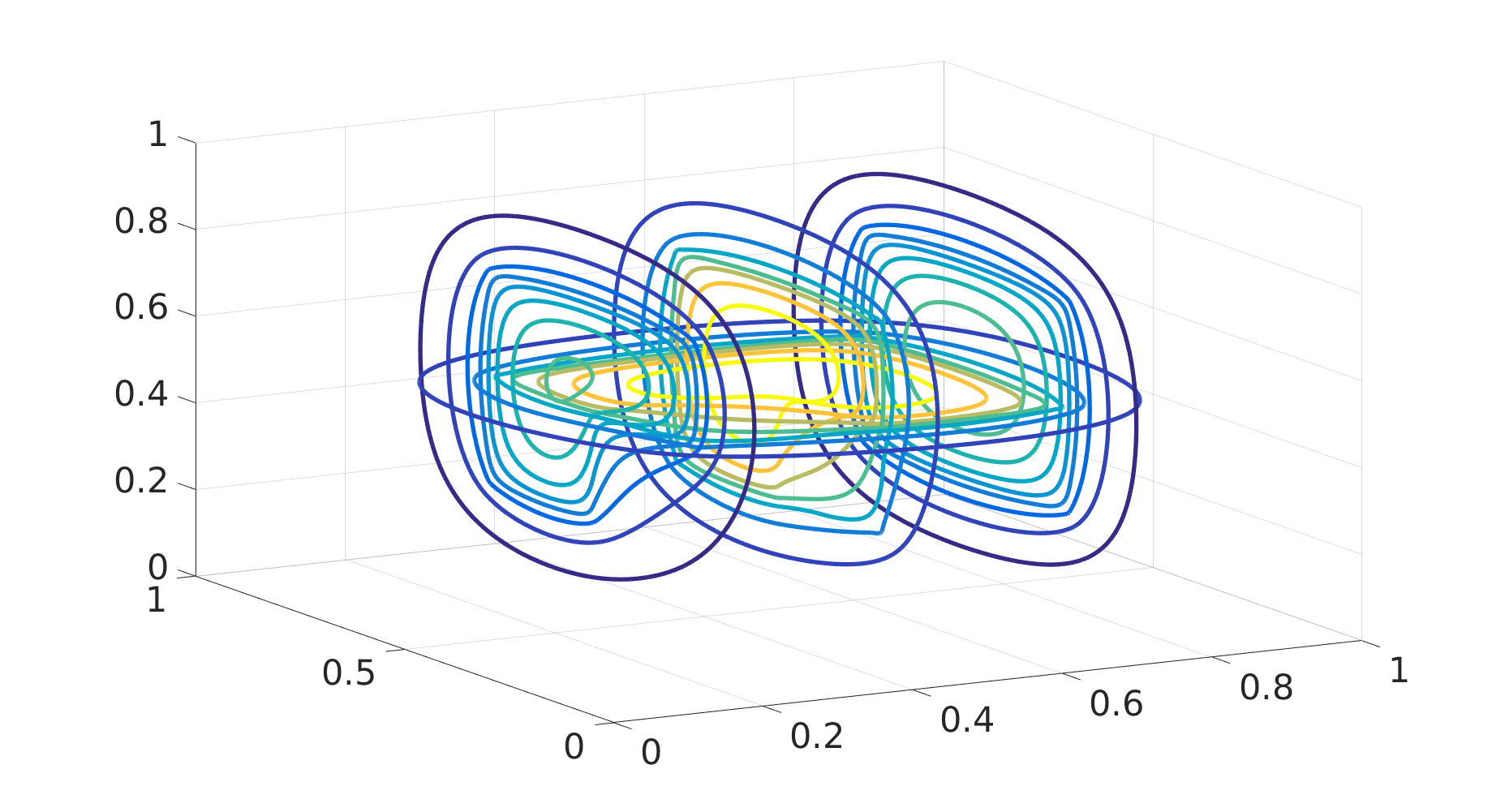}
\caption{Solutions of the equation with 3D right-hand sides (analogous to
Figure \ref{fig:design_RHS}) with $\alpha=1/2$ for $n=255$.}
\label{fig:FFT_3D_Alpha_12}
\end{figure}

  \begin{figure}[h]
\centering
\includegraphics[width=5.4cm]{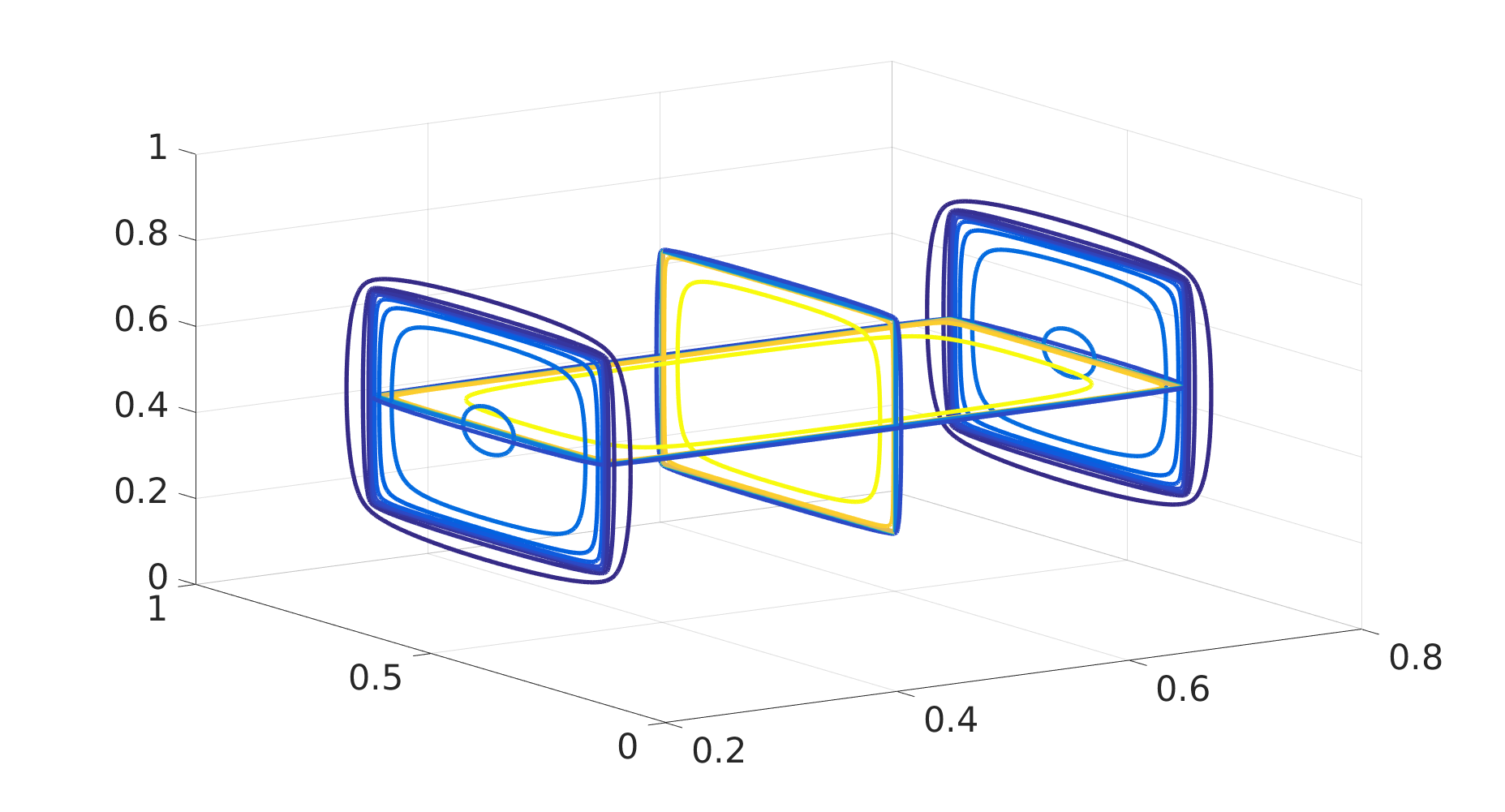}
\includegraphics[width=5.4cm]{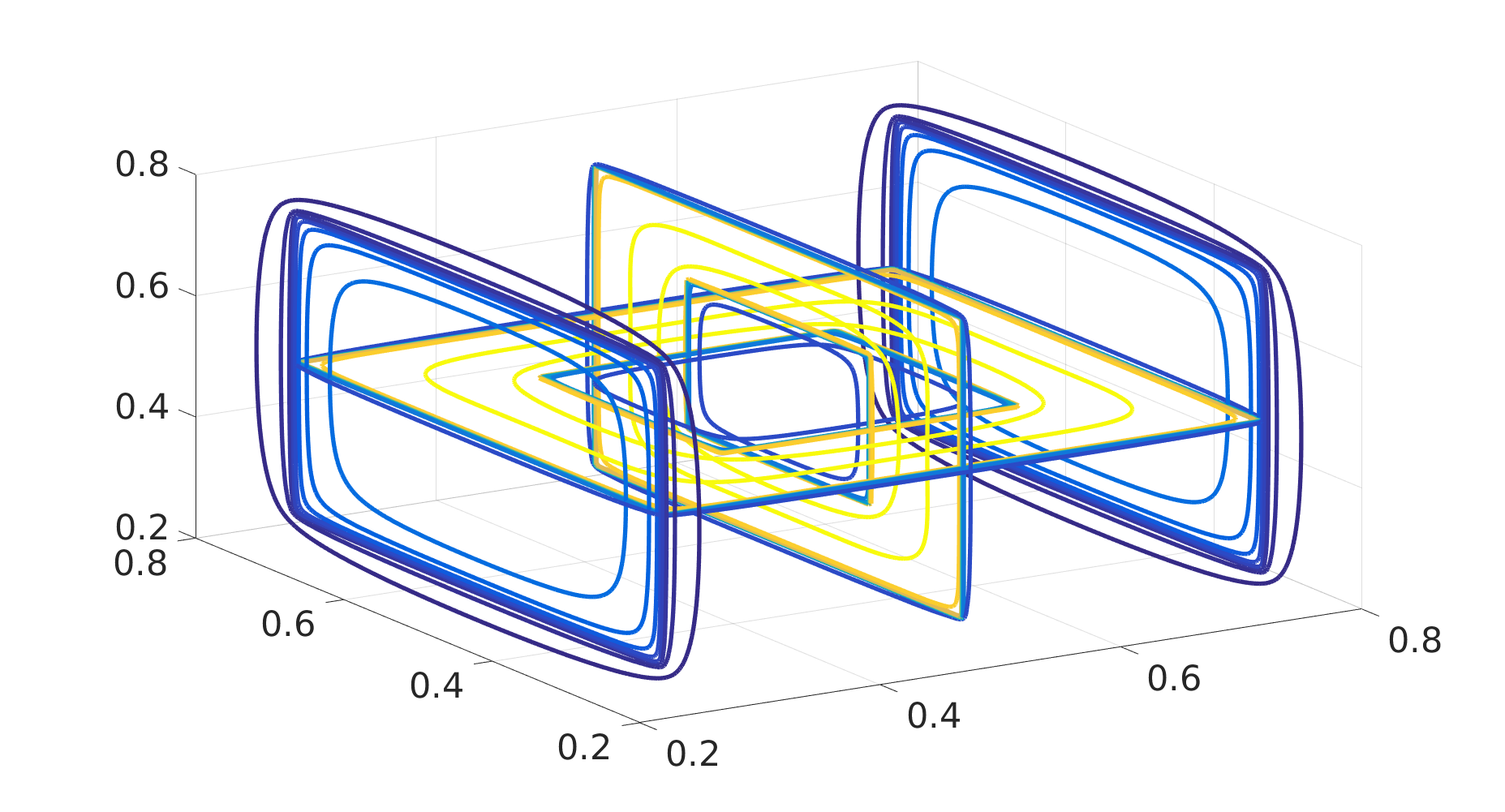}
\includegraphics[width=5.4cm]{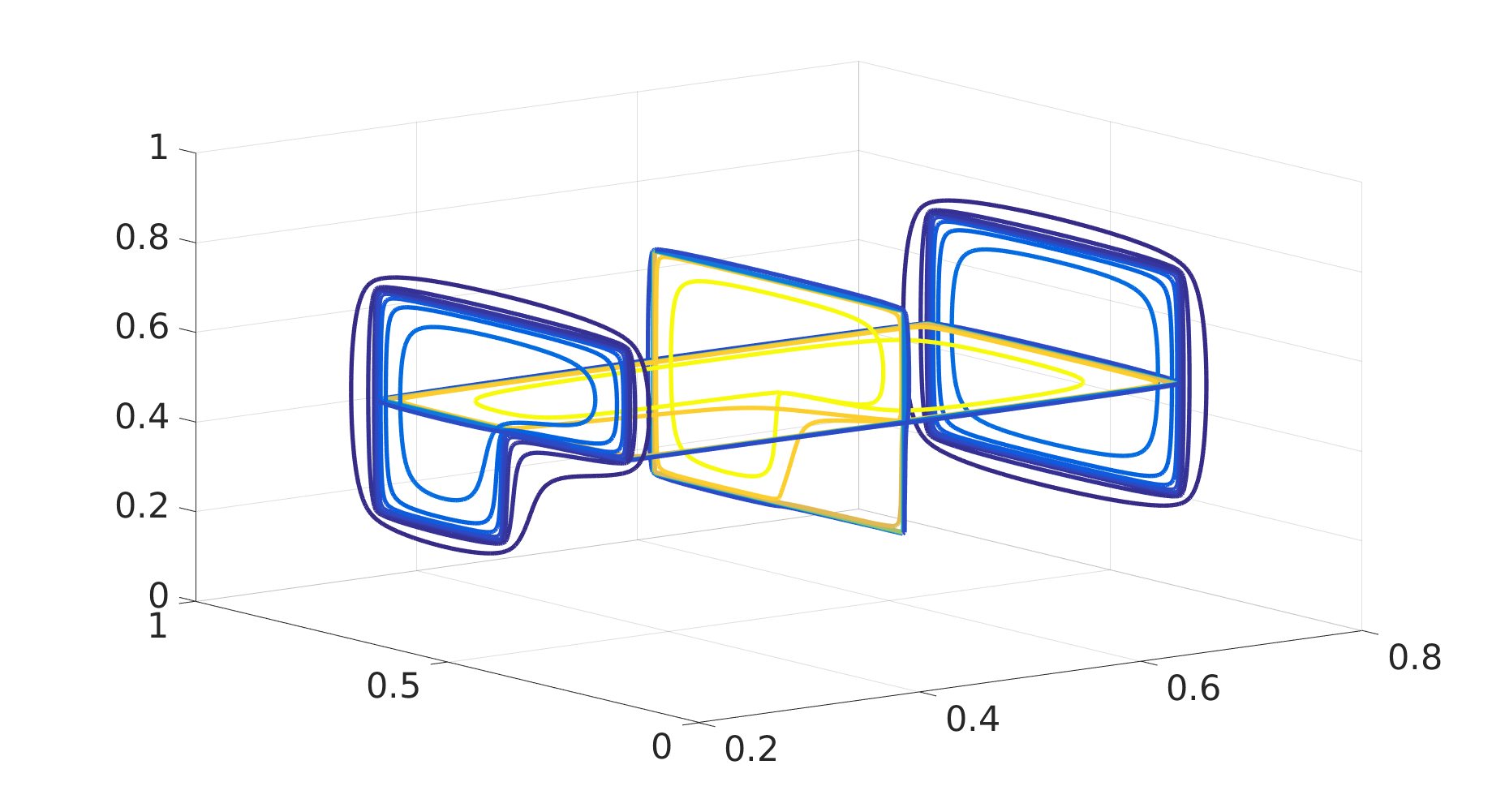}
\caption{\small Solutions of the equation with 3D right-hand sides (analogous to
Figure \ref{fig:design_RHS}) with $\alpha=1/10$ for $n=255$.}
\label{fig:FFT_3D_Alpha_110}
\end{figure}

\section{Conclusions}

We have introduced and analyzed a new approach for the optimal control of
a fractional Laplacian equation using tensor numerical methods.
The fractional Laplacian is diagonalized  in the FFT basis  on a tensor grid
and   a low Kronecker rank approximation to the core diagonal matrix is computed.
We present the novel rank-structured tensor approximation of functions of 
fractional elliptic operators based on sinc-approximation method.
 This representation exhibits the exponential decay of the
approximation error in the rank parameter.

These results apply to the fractional Laplacian itself, as well
as to the solution operators of a fractional control problem, resulting from
first-order necessary
conditions. Due to   the separation of spatial variables in tensor formats, the
application of the arising matrix-valued functions of the fractional
Laplacian to a given rank-structured vector has a complexity which is nearly linear
(linear-logarithmic) in the univariate grid size, independently of 
the spatial dimension of the problem.

The PCG iterative algorithm with adaptive rank truncation
for solving the equation for control function is implemented. In 3D case the rank
truncation is based on the RHOSVD-Tucker approximation and its transform to the low-rank
canonical form.
The numerical study illustrates the exponential decay of the
approximation error of the canonical tensor decompositions of the target tensors
in the rank parameter, 
and indicates the almost linear complexity scaling of the rank-truncated PCG solver
in the univariate grid size $n$ for 3D problems discretized on $n \times n \times n $
Cartesian grid. The PCG iteration exhibits the uniform
convergence rate in the univariate grid size $n$ and other model parameters.

 The tensor techniques presented in this paper can be generalized in several directions.
First of all, the approach can be extended to the elliptic operators with variable, 
but well separable, 
coefficients posed in the box type domains (e.g., layer type or perforated structures), 
which will be considered elsewhere. 
Further generalization to the case of non-rectangular domains is also possible. 
In this case one can use the alternative definition of the fractional elliptic operator 
by using the Dunford-Cauchy contour integral representation 
\cite{Higham_MatrFunc:08,HaHighTrefeth:08,GaHaKh3:02}
(see (\ref{eqn:Dunford_Int}) and related discussion) 
which is based on computations with the elliptic resolvent in a small number 
of quadrature points  located on the integration path.
The practical application of the Dunford-Cauchy representation requires the solution of linear systems 
of equations involving only the discrete elliptic operator (but not its fractional power). 
In this case the low-rank tensor decomposition techniques
can be applied on domains composed of the moderate number of box type subdomains 
(e.g., L-shaped, O-shaped or step-type domains).

Finally, we notice that the further reduction of the numerical complexity to the logarithmic scale
can be achieved by using the quantized-TT (QTT) representation of all discrete functions and 
operators involved, see \cite{KhQuant:09,Osel-TT-LOG:09,KhorBook:17}.

\begin{footnotesize}

\end{footnotesize}

\end{document}